\documentclass[12pt,onecolumn,oneside,leqno,a4paper]{article}
\usepackage{german}
\usepackage{exscale}
\usepackage{amsmath,amssymb,verbatim}
\usepackage{amscd}
\usepackage{graphicx}
\usepackage{rotating}
\usepackage{amsxtra,amsthm}
\usepackage{color}
\usepackage{makeidx}
\usepackage[pdftex]{hyperref}
\hypersetup{bookmarksopen=true,bookmarksopenlevel=1,bookmarksnumbered=true}

\newcommand{\TB}[1]{\textbf{#1}}
\newcommand{\TI}[1]{\textit{#1}}

\def\Bew{\textit{Proof. }}

\newcommand\QED{\hspace*{\fill} $ \square $}
\newcommand\euro{{\sffamily C%
  \makebox[0pt][l]{\kern-.70em\mbox{--}}%
  \makebox[0pt][l]{\kern-.68em\raisebox{.25ex}{--}}}}

\newenvironment{DES1}[1]
	{\begin{list}{}
	{\settowidth{\labelwidth}{#1}
	\setlength{\labelsep}{1em} 
	\setlength{\leftmargin}{0em}
	\setlength{\itemindent}{\labelwidth}
	\addtolength{\itemindent}{\labelsep}
	\setlength{\parsep}{\parskip}
	\setlength{\itemsep}{0.2\parsep plus0.02\parsep minus0.02\parsep}
	\setlength{\topsep}{0.2\parsep plus0.02\parsep minus0.02\parsep}
	\setlength{\partopsep}{0.6\parsep plus0.02\parsep minus0.05\parsep}
	}}
	{\end{list}}

\def\R{\mathbb{R}}

\def\Rplus{\mathbb{R}^+}
\def\Rplusnull{\mathbb{R}^+_0}

\def\B{\mathbb{B}}
\def\GL{\mathbb{GL}}

\def\PGL{\mathbb{PGL}}

\def\P{\mathbb{P}}

\def\SS{\mathbb{S}}
\renewcommand{\H}{\mathbb{H}}
\newcommand{\Rz}{\R^{2}}
\newcommand{\Rd}{\R^{3}}

\def\sign{\text{\textup{sign\,}}}

\def\id{\text{\textup{id}}}

\def\Span{\text{\textup{span}}}

\newcommand{\pfeil}[1]{\overrightarrow{#1}}

\newcommand\Forall{\;\forall\;}
\newcommand\Exists{\;\exists\;}
\newcommand\Impl{\;\Longrightarrow\;}

\newcommand{\so}{{\,}'} 

\newcommand\Arcsin{\arcsin\,}

\newcommand\Sin{\sin\,}
\newcommand\Cos{\cos\,}

\newcommand\vhi{\varphi}
\newcommand\vro{\varrho}
\newcommand\Sinh{\sinh\,}
\newcommand\Cosh{\cosh\,}
\newcommand\Tanh{\tanh\,}

\DeclareMathOperator{\arsinh}{arsinh}

\newcommand\Arsinh{\arsinh\,}

\renewcommand\Dot[2]{\left\langle #1,#2 \right\rangle}

\def\ds{\displaystyle}

\def\scr{\scriptscriptstyle}
\def\punkt{\scr \bullet}

\renewcommand{\Hat}[1]{\widehat{\,#1\,}}

\newcommand{\To}{\longrightarrow}

\newcommand{\cross}{\times}

\newcommand{\dun}{\begin{array}[b]{c} \punkt \\[-2ex] \cup \end{array}}
\newcommand{\Dunl}[2]
{\begin{array}[b]{c}
	{\scriptstyle \bullet}  \\[-3.8ex]
	\ds \bigcup_{#1}^{#2}
\end{array}}

\newcommand{\Sc}{\mathcal{S}}
\newcommand{\Dc}{\mathcal{D}}
\newcommand{\Hc}{\mathcal{H}}

\newcommand{\norm}[1]{\left\|#1\right\|}
\newcommand{\BM}{\begin{pmatrix}}
\newcommand{\EM}{\end{pmatrix}}
\newcommand{\BV}{\begin{vmatrix}}
\newcommand{\EV}{\end{vmatrix}}

\newtheoremstyle{tld}
{}
{}
{\itshape}
{}
{\bfseries}
{.}
{2em}
{}

\newtheoremstyle{bb}
{}
{}
{}
{}
{\bfseries}
{.}
{2em}
{}

\newtheoremstyle{bbs}
{}
{}
{}
{}
{\bfseries}
{.}
{\newline}
{}

\swapnumbers

\theoremstyle{tld}
\newtheorem{thm}{\negthickspace. Theorem}[section] 
\newtheorem{thmz}[thm]{\negthickspace. Zusatz}
\newtheorem{lem}[thm]{\negthickspace. Lemma}
\newtheorem{cor}[thm]{\negthickspace. Corollary}
\newtheorem{defi}[thm]{\negthickspace. Definition}
\newtheorem{defthm}[thm]{\negthickspace. Definition und Satz}
\newtheorem{thmdef}[thm]{\negthickspace. Satz und Definition}

\theoremstyle{bb}
\newtheorem{ex}[thm]{\negthickspace. Example}
\newtheorem{rem}[thm]{\negthickspace. Remark}
\newtheorem{aufg}{\negthickspace. Aufgabe}[section] 

\theoremstyle{bbs}
\newtheorem{exs}[thm]{\negthickspace. Examples}
\newtheorem{rems}[thm]{\negthickspace. Remarks}

\newcommand{\BT}{\begin{thm}}
\newcommand{\ET}{\end{thm}}

\newcommand{\BTZ}{\begin{thmz}}
\newcommand{\ETZ}{\end{thmz}}

\newcommand{\BL}{\begin{lem}}
\newcommand{\EL}{\end{lem}}

\newcommand{\BC}{\begin{cor}}
\newcommand{\EC}{\end{cor}}

\newcommand{\BD}{\begin{defi}}
\newcommand{\ED}{\end{defi}}

\newcommand{\BDT}{\begin{defthm}}
\newcommand{\EDT}{\end{defthm}}

\newcommand{\BTD}{\begin{thmdef}}
\newcommand{\ETD}{\end{thmdef}}

\newcommand{\BX}{\begin{ex}}
\newcommand{\EX}{\end{ex}}

\newcommand{\BXS}{\begin{exs}}
\newcommand{\EXS}{\end{exs}}

\newcommand{\BR}{\begin{rem}}
\newcommand{\ER}{\end{rem}}

\newcommand{\BRS}{\begin{rems}}
\newcommand{\ERS}{\end{rems}}

\newcommand{\BE}{\begin{equation}}
\newcommand{\EE}{\end{equation}}

\newcommand{\BA}{\begin{aufg}}
\newcommand{\EA}{\end{aufg}}

\numberwithin{equation}{section}

\newcommand{\bild}[1]{\includegraphics{BILDER/#1}}

\newcommand{\abs}[1]{\left|#1\right|}
\newcommand{\del}{\partial}

\begin{document} 

\thispagestyle{empty}

\begin{center}
\TB{\LARGE Polygons in hyperbolic geometry 1:}

\smallskip
\TB{\LARGE Rigidity and inversion of the $ n $-inequality}

\TB{\large Rolf Walter}
\end{center}

\setcounter{section}{0}


\section{\hspace{-1em}.
Introduction}
\label{intro}


In a recent series of papers, K. Leichwei"s extended various
topics from Euclidean geometry to the hyperbolic plane (and
partially to the geometry on the sphere). His results
concern convexity and extremum problems for closed curves,
including the addition of convex sets, support functions,
curves of constant width and Steiner's symmetrization; see
\TI{Leichtwei"s} [2003], [2004], [2005], [2008.a], [2008.b].

Looking in a different direction, we will discuss here
certain problems in the distance geometry of \TI{Menger}
[1928] and \TI{Blumenthal} [1970] for polygons in the
hyperbolic geometry. This part $ 1 $ deals with the rigidity
for polygons of fixed sidelengths and the 
inversion of the generalized triangle inequality. The
polygons here are of a general nature, no
a-priori assumptions on their form (convexity, simple
closedness, etc.) have to be made. For the Euclidean plane,
analogous problems were treated by \TI{Pinelis} [2005]. We
have been surprised to learn from the Pinelis paper that the
Euclidean solution was open for such a long time. Rigidity
is understood in the non-continuous sense similar to
\TI{Klingenberg} [1978], Theorem 6.2.8.

In the hyperbolic plane, as in Euclidean geometry, a polygon
usually can be changed in many ways while keeping its
sidelengths fixed. Exceptions are some special collinear
polygons and all polygons which are convex and cocyclic, the
latter meaning that the vertices lie on a circle. In fact,
we shall prove that these polygons can be changed in the
hyperbolic plane only in a trivial way, namely by rigid
hyperbolic motions if the sidelengths are kept constant
(Theorems \ref{thm_polyrig_haupt} and
\ref{thm_polyex_inv}). This rigidity problem
is more difficult than in the Euclidean case because there
are three different types of circles in the hyperbolic
plane, the distance circles, the distance lines, and the
horocycles. The solution depends on closedness conditions
for the length spectrum which are expressed in terms of the
generating groups of these figures.

A related question in the hyperbolic geometry is: Which
conditions on the length spectrum must be posed in order
that a corresponding polygon will exist? The answer is that
the generalized triangle inequality, henceforth called the 
$ n $-inequality, does the job (Corollary
\ref{cor_polyex_hn}). Again, the
closedness conditions are crucial for the proof.

Since both problems belong to distance geometry it is
necessary, to dispense to a certain extent with classical
hyperbolic trigonometry, because most of its relations are
between lengths \TI{and} angles. Expressing the situation
with lengths alone may render some formulation more
complicated, but provides a better insight into the
dependencies of quantities.

This aspect becomes even more important in part $ 2 $ where
the problem of maximizing the area of polygons with fixed
sidelengths will be treated, based on the results of the
present part $ 1 $. Moreover it will become clear that the
convex and cocyclic polygons are the only non-collinear ones
for which the rigidity can hold true.


\section{\hspace{-1em}.
The Cayley/Klein model}
\label{prelim}
\markright{\ref{prelim}. The Cayley/Klein model}

For our purposes, the \TI{Cayley/Klein model} is the
most suitable framework of hyperbolic geometry.
In this model, the \TI{points} of the hyperbolic plane are taken
as the ordinary points of the open unit ball
$$
\B := 
\{(\xi_{1},\xi_{2}) \in \R^2 \mid \xi_{1}^2+\xi_{2}^2 < 1\},
$$
and the \TI{straight lines} are taken as the ordinary
straight chords of $ \B $, always connecting two different
points of the \TI{horizon}
$$
\SS := \del \B =
\{(\xi_{1},\xi_{2}) \in \R^2 \mid \xi_{1}^2+\xi_{2}^2 = 1\},
$$
which itself does not belong to the hyperbolic plane.
Nevertheless, points of $ \R^2 $ outside $ \B $ may well be 
considered for describing situations within $ \B $. For
basics on hyperbolic geometry see e.g. \TI{Coxeter} [1968]
or \TI{Lenz} [1967].

The topology of $ \B $, its orientation and
between-relations on lines are the same as in standard 
$ \Rz $. The origin $ (0,0) \in \Rz $ is denoted by 
$ O $, the connecting line of two different points 
$ X,Y \in \Rz $ by $ X \vee Y $, and the directed line from 
$ X $ to $ Y $ by $ \pfeil{X \vee Y} $. The segment from
$ X $ to $ Y $ as a point set is denoted by 
$ \overline{XY} $; it includes $ X $ and $ Y $.
The open segment between $ X = (-1,0) $ and $ Y = (1,0) $ will be
called the \TI{groundline} and denoted by $ U_{0} $.

It is convenient to specify the 
points $ X = (\xi_{1},\xi_{2}) \in \Rd $ also by  
\TI{homogeneous coordinates} $ x_{0},x_{1},x_{2} $
defined by
$$
\xi_{1} = \frac{x_{1}}{x_{0}}, \qquad 
\xi_{2} = \frac{x_{2}}{x_{0}}, \qquad x_{0} \neq 0,
$$
thus imbedding $ \Rz $ in the projective plane $ \P^{2} $,
consisting of the one dimensional vector subspaces (rays) of 
$ \Rd $. 

All metric properties in $ \B $ are expressible by the
pseudo-Euclidean scalar product in $ \Rd $:
\BE
\label{eq_hcoo_iscal}
\Dot{x}{y} := x_{0}y_{0}-x_{1}y_{1}-x_{2}y_{2}, \qquad
x := (x_{0},x_{1},x_{2}), \quad y := (y_{0},y_{1},y_{2}).
\EE
For each $ x \in \Rd $ we use the abbreviation 
$$
\norm{x} := \sqrt{\abs{\Dot{x}{x}}},
$$
despite the fact that this does not define a norm in the
analytic sense.

In homogeneous coordinates, a line of $ \Rz $ has
an equation of the form $ \Dot{u}{x} = 0 $ 
($ u \neq 0 $ in $ \Rd $ fixed, 
$ x \neq 0 $ in $ \Rd $ variable). The points of $ \B $
are characterized by $ \Dot{x}{x} > 0 $ and the lines
hitting $ \B $ by $ \Dot{u}{u} < 0 $.

A vector $ x := (x_{0},x_{1},x_{2}) \neq 0 $ in $ \Rd $ 
belonging to a point of $ \Rz $ is called a 
\TI{(homogeneous) point vector}, a
vector $ u := (u_{0},u_{1},u_{2}) \neq 0 $ in $ \Rd $ 
belonging to a line in 
$ \Rz $ is called a \TI{(homogeneous) line vector}. Such a 
$ u $, viewed as a homogeneous point vector, represents the
polar point of the line w.r.t. to the circle $ \SS $. 
For a point vector $ x $ we always assume $ x_{0} > 0 $. If
convenient, one may $ x $ choose such that $ \Dot{x}{x} = 1 $
and $ u $ such that $ \Dot{u}{u} < -1 $. If so, we 
call $ x $ resp. $ u $ \TI{normalized}.
The normalized point vectors 
define a bijective correspondence of $ \B $ with the shell
of a hyperboloid in $ \Rd $
$$
\H := \{x \in \R^{3} \mid  \Dot{x}{x} = 1, \; x_{0} > 0\}.
$$
$ \H $ is the \TI{quadric model} of hyperbolic geometry. It 
may be considered parallel to the Cayley/Klein model 
what is sometimes helpful.

The scalar product \eqref{eq_hcoo_iscal} induces a
corresponding pseudo-Euclidean cross product from 
$ \Rd \cross \Rd $ to $ \Rd $ by the identity
\BE
\label{eq_hcoo_crpr}
[x,y,z] = \Dot{x \cross y}{z} \quad \Forall x,y,z \in \Rd,
\EE
where $ [x,y,z] $ denotes the standard determinant form on 
$ \Rd $. The cross product may serve to calculate a line 
vector $ u $ of the connecting line of two different points with 
point vectors $ x,y $
as $ u = x \cross y $ and a point vector $ x $ of the intersection
point of two different lines with line vectors $ u, v$ as 
$ x = u \cross v $.

Important rules, often to be used, are:
\begin{alignat}{2}
\label{eq_hcoo_gram}
[x,y,z] \cdot [x',y',z'] &=
\BV
\Dot{x}{x'} & \Dot{x}{y'} & \Dot{x}{z'} \\
\Dot{y}{x'} & \Dot{y}{y'} & \Dot{y}{z'} \\
\Dot{z}{x'} & \Dot{z}{y'} & \Dot{z}{z'}
\EV
&\quad&
\text{\TI{(Gram's identity)}}\\[1ex]
\label{eq_hcoo_lagr}
\Dot{x\cross y}{x'\cross y'} &=
\BV
\Dot{x}{x'} & \Dot{x}{y'} \\
\Dot{y}{x'} & \Dot{y}{y'}
\EV
&\quad&
\text{\TI{(Lagrange's identity)}}\\[1ex]
\label{eq_hcoo_tripl}
(x\cross y)\cross x' &= \Dot{x}{x'}y-\Dot{y}{x'}x
&\quad&
\text{\TI{(triple identity)}},
\end{alignat}
valid for all $ x,y,z,x',y',z' $ in the pseudo-Euclidean
space $ \Rd $.

\smallskip
If necessary, the various objects occurring so far may be
denoted more accurately in the following manner: points of 
$ \Rz $ as pairs of real numbers by $ X, Y,\ldots $ (also
as elements of the hyperbolic plane $ \B $),
corresponding point vectors by the corresponding small
letters $ x,y,\ldots $; straight lines of $ \Rz $ as point
sets by $ U,V,\ldots $ and corresponding line vectors by the
corresponding small letters $ u,v,\ldots $. Sometimes, 
by abuse of language, we are speaking of a point $ x $,
meaning the point $ X $ with point vector $ x $, and
similarly for straight lines.

The most important notions are that of the \TI{hyperbolic distance}
$ d(X,Y) $ of points $ X,Y \in \B $ and the \TI{hyperbolic angle}
$ \alpha(U,V) $ of two lines intersecting in $ \B $ by
\BE
\label{eq_prelim_metr}
\Cosh d(X,Y) = 
\frac{\Dot{x}{y}}{\norm{x}\cdot \norm{y}},
\qquad\qquad
\Cos \alpha(U,V) = 
\frac{\abs{\Dot{u}{v}}}%
{\norm{u}\cdot\norm{v}}, \quad
0 \leq \alpha \leq \frac{\pi}{2},
\EE
Again, this definitions only involve the pseudo-Euclidean
scalar product \eqref{eq_hcoo_iscal}. They are invariant
under the possible changes of the representing objects
$ x,y $, resp. $ u,v $. It is well known that $ d $ defines 
a complete metric on $ \B $.

If $ U \subset \B $ is a straight line and $ X \in \B $ any 
point then there is exactly one \TI{foot point} $ F \in U $,
i.e. $ F $ satisfying $ d(X,F) \leq d(X,Y) $ for all $ Y \in U $. 
The distance $ d(X,F) $ is the \TI{distance} $ d_{U}(X) $ 
of $ X $ and $ U $. Explicitly:
\BE
\label{eq_prelim_dU}
f := \Dot{u}{x}u-\Dot{u}{u}x,\qquad 
\Sinh d_{U}(X) = 
\frac{\abs{\Dot{u}{x}}}{\norm{u}\cdot\norm{x}}.
\EE
For $ X \notin U $ the foot $ F $ is also characterized by
the property that the connecting line of $ X $ and $ F $
cuts $ U $ orthogonally at $ F $.

The \TI{orientation} of $ \B $ shall be that inherited from the
standard orientation of $ \Rz $. The two 
\TI{half-spaces (sides)}
of a line $ U $ are defined by $ \Dot{u}{z} > 0 $,
resp. $ \Dot{u}{z} > 0 $. If $ u $ is multiplied by a negative
factor then both sides are interchanged, so no side is
distinguished. However, if $ \pfeil{X \vee Y} $ is the
directed line from $ X $ to $ Y $ (where $ X \neq Y $) then 
the line vector $ u := x \cross y $ is fixed up to a 
\TI{positive} factor. Thus a directed line may also be given 
by a line vector $ u $ and its positive multiples.
In this case the two sides of $ U $ can be distinguished by 
$ \Dot{u}{z} > 0 $, resp. $ \Dot{u}{z} < 0 $. The first one is
called the \TI{left-hand side}, the second one the 
\TI{right-hand side} of the directed line. For 
$ u = x \cross y $ these inequalities are equivalent to 
$ [x,y,z] > 0 $, resp. $ [x,y,z] < 0 $.

\medskip
\TB{\TI{Hyperbolic isometries}}

The isometries of the hyperbolic plane $ \B $ can be
generated from a rational parametrization of the horizon $ \SS $
by re-parametrizing $ \SS $ with a real broken linear transformation 
(a \TI{homography} of $ \R $) and extending this
transformation of $ \SS $ projectively to $ \B $.
This amounts to the following: Associate to any real 
$ (2 \cross 2) $-matrix $ A $ a real $ (3 \cross 3) $-matrix 
$ M(A) $ by
\BE
\label{eq_hypmot_M}
A:=
\BM
a & b \\
c & d
\EM
\longmapsto
M(A):=
\BM
a^{2}+b^{2}+c^2+d^{2} & -a^{2}+b^{2}-c^2+d^{2} &
2(ab+cd)\\
-a^{2}-b^{2}+c^2+d^{2}&a^{2}-b^{2}-c^2+d^{2}&
2(-ab+cd)\\
2(ac+bd)&2(-ac+bd)&2(ad+bc)
\EM.
\EE
Then simple (though tedious) matrix calculations show the
following rules for all real $ (2 \cross 2) $-matrices 
$ A,B $ and all vectors $ x,y,z \in \Rd $:
\begin{align}
\label{eq_gt_hom1}
M(A\cdot B) &= 2\cdot M(A) \cdot M(B) \\[1ex]
\label{eq_gt_det}
\det M(A) &= 8(\det A)^{3} \\[1ex]
\label{eq_gt_detb}
[M(A)x,M(A)y,M(A)z] &= 8(\det A)^{3}[x,y,z] \\[1ex]
\label{eq_gt_invsk1}
\Dot{M(A)x}{M(A)y} &= 4(\det A)^{2}\Dot{x}{y} \\[1ex]
\label{eq_gt_invskc}
(M(A)x)\cross(M(A)y) &= 2\det A \cdot  M(A)(x \cross y).
\end{align}
In particular, if $ \det A = ad-bc \neq 0 $, then $ M(A) $
defines a non-singular linear map from $ \Rd $ onto itself
and thus a bijective projective map from $ \P^{2} $
to itself. In addition, the rule \eqref{eq_gt_invsk1} signifies that
this projective map
induces an isometry $ \mu(A) $ of the hyperbolic plane 
$ \B $ onto itself. Each such transformation $ \mu(A) $ will be 
called a \TI{hyperbolic isometry} 
of $ \B $.
It transforms lines to lines and leaves invariant the 
hyperbolic metric.
The $ \mu(A) $ is on $ \B $ orientation-preserving if 
$ \det(A) > 0 $ (equivalently: $ \det M(A) > 0 $) and
orientation-reversing if $ \det(A) < 0 $ (equivalently: 
$ \det M(A) < 0 $).

Moreover, the rule \eqref{eq_gt_hom1} shows that the
isometries of $ \B $ form a group homomorphic via $ \mu $ to the 
linear group $ \GL(2,\R) $ modulo its center, 
i.e. to the real projective
group $ \PGL(2,\R) $ of the projective extension of $ \R $. 
In fact, this 
homomorphism is an isomorphism because $ M(A) $ being
a multiple of the unit $ (3 \cross 3) $-matrix implies
that $ A $ is a multiple of the unit $ (2 \cross 2) $-matrix.

An orientation-preserving isometry $ \mu(A) $ will be called 
a hyperbolic \TI{motion}. The motions form a subgroup of
index $ 2 $ of the group of all isometries. In order to fix 
the notation for the present purposes, two figures in $ \B $
will be called \TI{hyperbolically equivalent}, resp. \TI{congruent} 
if they are transformed one to each other by a hyperbolic
motion, resp. a hyperbolic isometry. In most cases the
claims will concern equivalence, i.e. the stronger version.
Sometimes, the adjective `hyperbolic' is skipped.

\medskip
\TB{\TI{Hyperbolic circles}}

The straight lines in $ \B $ correspond to one-dimensional
sections of the quadric shell $ \H $ with vector planes
in $ \Rd $. Analogous to sphere geometry, circles in the
hyperbolic plane $ \B $ are defined by correspondence with
the possible intersections of $ \H $ with planes of 
$ \Rd $, not containing the origin, if the intersections
span the planes. Such a plane has an equation of the form 
$ \Dot{u}{x} = p $ with fixed $ u \neq 0 $ in $ \R^{3} $
and $ p \neq 0 $ in $ \R $ and additional conditions to 
ensure that the intersection is really a conic section.
This induces a classification as follows: 

In case $ \Dot{u}{u} > 0 $, the quadratic form 
$ \Dot{\;}{\;} $ is negative definite on the orthogonal
complement $ u^{\bot} $ of $ \Span(u) $, and the
intersection is an ellipse if it spans the plane. Without
loss, one may assume $ u_{0} > 0 $. The ellipse exists and
is contained in the half-space $ x_{0} > 0 $ of $ \Rd $ iff 
$ p > \norm{u} $. Then the pre-image in $ \B $ of the intersection  
is called a \TI{distance circle}.

In case $ \Dot{u}{u} < 0 $, the quadratic form 
$ \Dot{\;}{\;} $ is indefinite on $ u^{\bot} $, and the 
intersection always is a hyperbola with one component in the
half-space $ x_{0} > 0 $. Then the pre-image in $ \B $ of the intersection  
is called a \TI{distance line}.

In case $ \Dot{u}{u} = 0 $, the orthogonal space $ u^{\bot} $
does not complete $ \Span(u) $ (since 
$ u \in u^{\bot} $). The intersection always is a parabola. 
Again it can be assumed $ u_{0} > 0 $. The parabola is
contained in the half-space $ x_{0} > 0 $ iff $ p > 0 $. 
Then the pre-image in $ \B $ of the intersection 
is called a \TI{horocycle}.

Each of these pre-images in $ \B $ is called a \TI{circle}.
From this analytic definition  the geometric
meanings of circles are easily deduced as follows:

\medskip
\BL[data of the hyperbolic circles]
\label{lem_circ_A}
~\\[0.3ex]
From the circle equation $ \Dot{u}{x} = p $ of the
foregoing classification, the non-Euclidean data of
circles can be read off:
\begin{DES1}{(iii)}
\item[(i)]
For a distance circle ($ \Dot{u}{u} > 0 $) let $ c $ and 
$ R > 0 $ be defined by
\BE
\label{eq_circ_A1}
c := \frac{u}{\norm{u}}, \qquad \Cosh R := \frac{p}{\norm{u}}.
\EE
Then the distance circle consists of all points $ X \in \B $ 
with $ d(C,X) = R $. So  $ C \in \B $ is the \underline{center}
and $ R $ is the \underline{(hyperbolic) radius}.
\item[(ii)]
For a distance line ($ \Dot{u}{u} < 0 $), let $ \delta > 0 $
be defined by
\BE
\label{eq_circ_A2}
\Sinh \delta := \frac{\abs{p}}{\norm{u}}.
\EE
Then the distance line consists of all points 
$ X \in \B $ with $ d_{U}(X) = \delta $, situated on the
side of $ U $ with $ \Dot{u}{x} = \sign p $. The line 
$ U \subset \B $ will be called the \underline{soul} of the distance
line. 
\item[(iii)]
For a horocycle ($ \Dot{u}{u} = 0 $), the vector $ u $, if
considered as a point vector, defines a point $ U $ on the
horizon. The horocycle is an orthogonal trajectory of all
straight lines of $ \B $ passing through $ U $. 
\QED
\end{DES1}
\EL

\smallskip
\BR
\label{rem_circ_d23}
The hyperbolic circles in $ \B $, when viewed with
`Euclidean eyes', are elliptic arcs with the following
design in standard positions:

\begin{center}
\begin{minipage}{15cm}
\bild{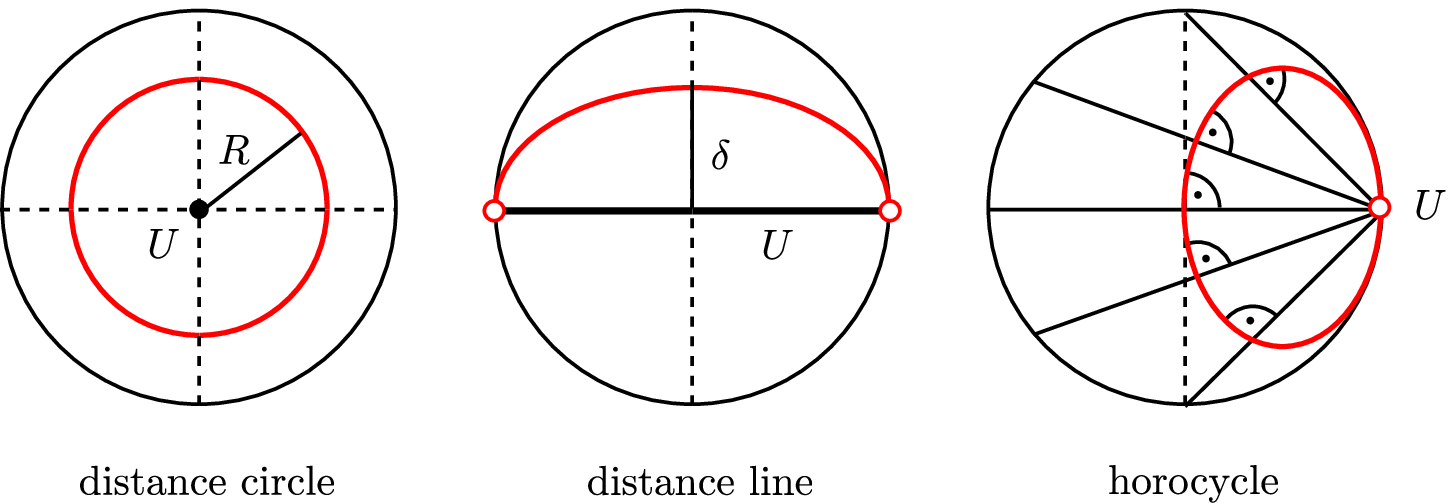}
\end{minipage}
\end{center}

A distance circle with center $ C = (0,0) $ is an ordinary 
circle with this center. The relation between the hyperbolic
radius $ R $ and the Euclidean radius $ r $ is:
$$
r = \Tanh R.
$$
A distance line whose soul is the groundline $ U_{0} $ is a
halfellipse with midpoint $ O $ and axes parallel to the
coordinate axes with Euclidean half-axes $ a = 1 $ resp. 
$ b = \Tanh \delta $, and situated either in the upper or
lower halfplane. A horocycle with horizon point 
$ U = (1,0) $ is an elliptic arc from $ U $ to $ U $ again
having its axes parallel to the coordinate axes where the
half-axes are in the relation $ a = b^{2} $, e.g. 
$ a = \frac{1}{2} $, $ b = \frac{\sqrt{2}}{2}$.
\ER

It will become essential in our investigation that the various circles are 
orbits of one parameter groups of hyperbolic motions. These 
groups arise as follows: One starts with the one
parameter groups of the projective group $ \PGL(2,\R) $ via 
the Jordan normal forms of their infinitesimal generators.
Then with the homomorphism \eqref{eq_hypmot_M} one translates these 
groups into the hyperbolic plane. As infinitesimal
generators in the group $ \PGL(2,\R) $ one may take those in 
the group $ \GL(2,\R) $ with trace $ 0 $. These matrices 
can even be multiplied by a scalar factor because this only 
influences the parametrization. So, finally, three normal
forms for the infinitesimal generators survive, namely of the
types `complex diagonalizable', `real diagonalizable', and 
`nilpotent': 
\BE
\label{eq_hypmot_JNF}
\BM
0 & -1 \\
1 & \phantom{-}0
\EM,
\qquad
\BM
1 & 0 \\
0 & -1
\EM,
\qquad
\BM
0 & 0 \\
1 & 0
\EM.
\EE
The corresponding groups themselves arise from these
matrices by multiplication with the group parameter $ t $ and
exponentiation as:
$$
A(t): \quad
\BM
\Cos t & -\Sin  t \\
\Sin  t & \phantom{-}\Cos  t 
\EM,
\qquad
\BM
e^{t}& 0 \\
0 & e^{-t}
\EM,
\qquad
\BM
1 & 0 \\
t & 1
\EM.
$$
Inserting this in the homomorphism \eqref{eq_hypmot_M} then 
yields (after suitable abbreviations and re-scaling of 
$ M(A) $):
\BE
\label{eq_hypmot_1pg}
\BM
1 & 0 & 0 \\[1ex]
0 & \Cos \vhi & -\Sin \vhi \\[1ex]
0 & \Sin \vhi & \phantom{-}\Cos \vhi
\EM,
\qquad
\BM
\cosh \vro & \sinh \vro & 0 \\[1ex]
\sinh \vro & \cosh \vro & 0 \\[1ex]
0 & 0 & 1
\EM,
\qquad
\BM
1+\frac{1}{2}b^2 & -\frac{1}{2}b^2 & -b \\[1ex]
\frac{1}{2}b^2 & 1-\frac{1}{2}b^2 & -b \\[1ex]
-b & b & 1
\EM.
\EE
The new variables $ \vhi $, $ \vro $, $ b $ replacing 
$ t $ are also `good' group parameters, i.e. the composition in 
the group corresponds to the addition of the parameters.
Eqns. \eqref{eq_hypmot_1pg} are the desired normal forms
for the one parameter groups of motions, represented as 
linear groups in $ \Rd $, respecting the pseudo-Euclidean
scalar product up to factors. The groups induced on $ \B $
by \eqref{eq_hypmot_1pg} are known resp. as 
\TI{hyperbolic rotations} around $ O $, 
\TI{hyperbolic translations} along the groundline $ U_{0} $, and
\TI{limit rotations} with horizon point $ (1,0) \in \SS $.
A hyperbolic rotation around $ O $ with the parameter 
$ \vhi $ is the same as the Euclidean rotation around $ O $
with angle $ \vhi $. A hyperbolic translations along the 
groundline $ U_{0} $ induces on $ U_{0} $ a bijection with
the property that the distance between pre-image and image 
is constant.

Now, if $ (S(t))_{t\in\R} $ is such a representation (with
neutral notation $ t $ for the group parameter) one has the 
differential equation $ S' = GS = SG $ where $ G $ 
is the infinitesimal generator $ G := S'(0) $. The orbits
are described by $ \gamma(t) := S(t)a $ with $ a $ fixed in 
$ \Rd $. An orbit $ \gamma $ is contained in the quadric shell 
$ \H $ iff $ a \in \H $. For the derivatives of $ \gamma $ 
one has $ \gamma^{(k)} = G^{k}\gamma $ and also 
$ \gamma^{(k)} = SG^{k}a $. This implies 
$ \Span(\gamma\so,\gamma\so',\ldots,\gamma^{(k)}) = 
S(\Span(Ga,G^{2}a,\ldots,G^{k}a)) $, so the
dimensions of all osculating spaces remain constant along
the orbit. If especially $ Ga, G^{2}a $ are linearly
independent but $ Ga, G^{2}a, G^{3}a $ are linearly
dependent, then the orbit is contained in an affine plane
of $ \Rd $. The corresponding vector $ u $ of this plane is 
given by $ u := Ga \cross G^{2}a $ and, if 
$ \Dot{u}{a} \neq 0 $, the orbit is part of a circle, whose 
type is determined by the sign of
\BE
\label{eq_hypmot_sign}
\Dot{u}{u} = 
\Dot{Ga}{Ga}\Dot{G^{2}a}{G^{2}a}-\Dot{Ga}{G^{2}a}^{2}.
\EE
The condition $ \Dot{u}{a} \neq 0 $ is equivalent to 
$ [a,Ga,G^{2}a] \neq 0 $ (what already implies the linear 
independency of $ Ga,G^{2}a $). The foregoing assumption on 
linear dependency is satisfied if $ G $ is singular
because $ [Ga,G^{2}a,G^{3}a] = \det(G)\cdot [a,Ga,G^{2}a] $.
For the one parameter groups \eqref{eq_hypmot_1pg} this is fulfilled
since their infinitesimal generators are
\BE
\label{eq_hypmot_erz}
G: \qquad
\BM
0 & 0 & 0 \\
0 & 0 & -1 \\
0 & 1 & 0
\EM,
\qquad
\BM
0 & 1 & 0 \\
1 & 0 & 0 \\
0 & 0 & 0
\EM,
\qquad
\BM
0 & 0 & -1 \\
0 & 0 & -1 \\
-1 & 1 & 0
\EM.
\EE
This discussion shows that for our three one parameter
groups the resulting orbits solely depend on the choice of
the initial point $ a $. However, for initial points 
equivalent under an element of the group, the whole 
orbits will be equivalent in the same manner. The
classification of the orbits is thus reduced further.
Before formulating the result, the arc lengths of 
the orbits may be calculated.

\medskip
The pseudo-Euclidean scalar product in $ \R^{3} $ induces a 
corresponding pseudo-riemannian metric on $ \R^{3} $:
\BE
\label{eq_hypmot_metH} 
\Dot{dx}{dx} = (dx_{0})^2-(dx_{1})^2-(dx_{2})^2.
\EE
This metric becomes strictly negative definite when
restricted to the quadric shell $ \H $. Denoting this
restriction with the same symbol, a genuine riemannian metric
on $ \H $ is then $ -\Dot{dx}{dx} $. Pulling back this metric 
onto $ \B $ yields the riemannian metric of the hyperbolic
plane. In the cartesian coordinates of $ \B $ the pull back 
is calculated as 
$$ 
ds^{2} = 
\frac{1}{(1-\xi_{1}^{2}-\xi_{2}^{2})^{2}}
\Bigl(
(1-\xi_{2}^{2})d\xi_{1}^{2}+2\xi_{1}\xi_{2} d\xi_{1}
d\xi_{2}+(1-\xi_{1}^{2})d\xi_{2}^{2}
\Bigr).
$$
Now, if as above $ \gamma = Sa $ is an orbit on $ \H $ of a 
one parameter group $ (S(t))_{t\in \R} $ one has
$$
-\Dot{\gamma\so}{\gamma\so} = -\Dot{S'a}{S'a} = 
-\Dot{SCa}{SCa} = -\Dot{Ca}{Ca}.
$$
From this follows that the orbit is automatically
parametrized proportional to arclength, the arclength
integral being
\BE
\label{eq_hypmot_arcl}
\int_{t_{1}}^{t_{2}} \sqrt{-\Dot{Ca}{Ca}}\;dt = 
\norm{Ca}(t_{2}-t_{1}).
\EE
Based on all this the following statements can be proved.
The main concern is the relation between arclengths and
chordlengths which will become crucial in the discussion of
the closedness conditions in Sect. \ref{closed}.

\medskip
\BL
\label{lem_hypmot_LS}
~\\[-3ex]
\begin{DES1}{(iii)}
\item[(i)]
The distance circle $ \Sc(R) $ with center $ O $ and
hyperbolic radius $ R $ is identical with the image of the
orbit of the first group in \eqref{eq_hypmot_1pg} with
initial point $ a := (\Cosh R,\Sinh R,0) $. This orbit is
represented by
\BE
\label{eq_hypmot_LS1}
\gamma(\vhi) := 
\BM
\Cosh R \\
\Sinh R\cdot\Cos \vhi \\
\Sinh R\cdot\Sin \vhi
\EM.
\EE
When restricted to any halfopen interval $ P $ of length 
$ 2\pi $, this orbit provides a bijective map of $ P $ onto 
$ \Sc(R) $. Define the function 
\BE
\label{eq_hypmot_LS2}
\Phi(R,t):=
2\Sinh R\cdot\arcsin
\left(\frac{\ds\Sinh \frac{t}{2}}{\Sinh R}\right),
\qquad R > 0, \quad 0 \leq t \leq 2R,
\EE
Then, for any two points $ X \neq Y $ in $ \Sc(R) $, the
following relation is valid between the chordlength $ L := d(X,Y) $ 
and the arclength $ s $ of $ \Sc(R) $ along the arc of 
$ \Sc(R) $ residing on the right hand side of the directed
line $ \pfeil{X \vee Y} $:
\BE
\label{eq_hypmot_LS3}
s =
\begin{cases}
\Phi(R,L) & 
\text{if $ O $ is left of or on $ \pfeil{X\vee Y} $} \\[1ex]
2\pi\Sinh R-\Phi(R,L) & 
\text{if $ O $ is right of $ \pfeil{X\vee Y} $}.
\end{cases}
\EE
\item[(ii)]
The distance line $ \Dc(\delta) $ of distance $ \delta > 0 $
from the groundline situated in the upper halfplane is
identical with the image of the orbit of the second group 
in \eqref{eq_hypmot_1pg} with initial point
$ a := (\Cosh \delta,0,\Sinh \delta) $. This orbit 
\BE
\label{eq_hypmot_LS4}
\gamma(\vro) := 
\BM
\Cosh \delta\cdot\Cosh\vro \\
\Cosh \delta\cdot\Sinh\vro \\
\Sinh \delta
\EM, \qquad \vro \in \R,
\EE
provides a bijective map from $ \R $ onto $ \Dc(\delta) $.
Define the function 
\BE
\label{eq_hypmot_LS5}
\Psi(\delta,t):=
2\Cosh \delta\cdot\arsinh
\left(\frac{\ds\Sinh \frac{t}{2}}{\Cosh \delta}\right),
\qquad \delta > 0, \quad t \geq 0.
\EE
Then, for any two points $ X \neq Y $ in $ \Dc(\delta) $ the
following relation is valid between the chordlength 
$ L := d(X,Y) $ and the arclength $ s $ on $ \Dc(\delta) $
from $ X $ to $ Y $:
\BE
\label{eq_hypmot_LS6}
s = \Psi(\delta,L)
\EE
\item[(iii)]
The horocycle $ \Hc $ with horizon point $ (1,0) \in \Rz $
through the origin $ O $ is 
identical with the image of the orbit of the third group 
in \eqref{eq_hypmot_1pg} with initial point
$ a := (1,0,0) $. This orbit
\BE
\label{eq_hypmot_LS7}
\gamma(b) := 
\BM
1+\frac{1}{2} b^{2}\\[1ex]
\frac{1}{2} b^{2} \\[1ex]
-b
\EM, \qquad b \in \R,
\EE
provides a bijective map from $ \R $ onto $ \Hc $.
Define the function 
\BE
\label{eq_hypmot_LS8}
\Omega(t) := 2\sinh \frac{t}{2}, \qquad t \geq 0.
\EE
Then, for any two points $ X \neq Y $ in $ \Hc $ the
following relation is valid between the chordlength 
$ L := d(X,Y) $ and the arclength $ s $ on $ \Hc $
from $ X $ to $ Y $:
\BE
\label{eq_hypmot_LS9}
s = \Omega(L).
\EE
The horocycles don't have invariants, i.e. all horocycles
are hyperbolically equivalent to $ \Hc $.
\end{DES1}
\EL

\Bew
The orbits of the three groups \eqref{eq_hypmot_1pg} are
located in planes of $ \Rd $. In order to
identify the orbits with the various circles one only has 
to specify the initial point  $ a $ by the regulations from 
above 
($ \Dot{u}{a} = [a,Ca,C^{2}a] \neq 0 $ with $ C $ from
\eqref{eq_hypmot_erz}, test on the sign of $ \Dot{u}{u} $ 
according to \eqref{eq_hypmot_sign}, reduction
of $ a $ modulo the group action). In detail:

\smallskip
1) \TI{For the first group \eqref{eq_hypmot_1pg}:}

For an arbitrary initial point $ a \in \H $ one calculates 
$ \Dot{u}{a} = a_{0}(a_{1}^{2}+a_{1}^{2}) $ and
$ \Dot{u}{u} = (a_{1}^{2}+a_{1}^{2})^2 $. The reduction of 
$ a $ is then possible to $ a_{2} = 0 $, i.e.
$$
a := 
\frac{1}{\sqrt{1-r^2}}
\BM
1\\r\\0
\EM, \quad
0 < r < 1,
\quad\text{with the orbit}\quad
\gamma(\vhi) := 
\frac{1}{\sqrt{1-r^2}}
\BM
1\\r\Cos \vhi\\r\Sin \vhi
\EM
$$
where all conditions on $ a $ are fulfilled. The distance 
of the orbit points in $ \B $ from $ O $ is constant and
given by $ \Cosh d(O,\Gamma(\vhi)) =  (1-r^{2})^{-1/2} $.
The distance $ d(O,\Gamma(\vhi)) $ becomes equal to a given 
number $ R > 0 $ iff $ r = \Tanh R $. Thus, the various
orbits can attain every hyperbolic radius $ R $. 

For $ R $ fixed, the representation of $ \gamma(\vhi) $ shows
that the orbit points cover the whole distance circle 
$ \Sc(R) $ with center $ O $ and hyperbolic radius $ R $,
and they reach every point exactly once if $ \vhi $
is restricted on a halfopen interval of length $ 2\pi $.

From \eqref{eq_hypmot_arcl}, the 
arclength on $ \Sc(R) $ between two points 
$ \Gamma(\vhi_{1}) $ and $ \Gamma(\vhi_{2}) $ is calculated to be
\BE
\label{eq_hypmot_sf}
\sigma = \Sinh R \cdot \abs{\vhi_{1}-\vhi_{2}}.
\EE
In particular the perimeter of $ \Sc(R) $ is $ 2\pi \Sinh R $.
Further, the chordlength $ L $ between two arbitrary points 
$ \Gamma(\vhi_{1}) $ and $ \Gamma(\vhi_{2}) $ is given by 
\BE
\label{eq_hypmot_Lf}
\Sinh \frac{L}{2} = 
\Sinh R \cdot \Sin \frac{\abs{\vhi_{1}-\vhi_{2}}}{2}.
\EE
For points $ X \neq Y $ on $ \Sc(R) $ in a given order, the
arc $ B(X,Y) $ associated to the \TI{directed} chord from 
$ X $ to $ Y $ shall consist of all points of $ \Sc(R) $ on the
right-hand side of the directed line $ \pfeil{X \vee Y} $.
Now let $ s $ always be the hyperbolic arclength measure of 
$ B(X,Y) $. By \eqref{eq_hypmot_sf}, for $ X = \Gamma(\vhi_{1}) $ 
and $ Y = \Gamma(\vhi_{2}) $, one has to calculate $ s $ 
as $ s = \sigma $ if $ O $ is left of or on 
$ \pfeil{X \vee Y} $, i.e. if $ \vhi_{1}, \vhi_{2} $ are
chosen with 
$ \vhi_{1} < \vhi_{2} $ and $ \vhi_{2}-\vhi_{1} \leq \pi $,
but as $ s =  2\pi \Sinh R - \sigma $ if $ O $ is right of 
$ \pfeil{X \vee Y} $, i.e. if $ \vhi_{1}, \vhi_{2} $ are 
chosen with 
$ \pi < \vhi_{2}-\vhi_{1} < 2\pi $. Eliminating in both
cases the difference $ \abs{\vhi_{1}-\vhi_{2}} $ from the
equations for $ s $ and $ L $ yields the assertion 
\eqref{eq_hypmot_LS3}.

\smallskip
2) \TI{For the second group \eqref{eq_hypmot_1pg}:}

For an arbitrary initial point $ a \in \H $ one calculates:
$ \Dot{u}{a} = a_{2}(a_{1}^{2}-a_{0}^{2}) $, 
$ \Dot{u}{u} = -(a_{1}^{2}-a_{0}^{2})^2 $. The reduction of 
$ a $ is then possible to reach $ a_{1} = 0 $, i.e.
$$
a := 
\frac{1}{\sqrt{1-\eta^2}}
\BM
1\\0\\ \eta
\EM, \quad
0 < \abs{\eta} < 1,
\quad\text{with the orbit}\quad
\gamma(\vro) := 
\frac{1}{\sqrt{1-\eta^2}}
\BM
\Cosh \vro \\ \Sinh \vro \\ \eta
\EM
$$
where all conditions on $ a $ are satisfied with 
$ u = (0,0,1) $, the line vector of the groundline $ U_{0} $.

The distance of the orbit points in $ \B $ from $ U $
is independent of $ \vro $ and given by 
$ \Sinh d_{U}(\Gamma(\vro)) =
\abs{\eta}(1-\eta^{2})^{-1/2} $. 
The distance $  d_{U}(\Gamma(\vro)) $ becomes equal to a given 
number $ \delta > 0 $ iff $ \abs{\eta} = \Tanh \delta $. Thus, 
the various orbits can attain every hyperbolic distance 
$ \delta $ from $ U $.

For $ \delta $ fixed, the representation of $ \gamma(\vro) $ shows
that the orbit points run through the whole distance line in the
upper resp. lower halfplane (according to the sign of 
$ \eta $) reaching each point exactly once. Let 
$ \Dc(\delta) $ be the upper distance line. 

From \eqref{eq_hypmot_arcl}, the 
arclength on $ \Dc(\delta) $ between two points 
$ \Gamma(\vro_{1}) $ and $ \Gamma(\vro_{2}) $ is calculated to be
\BE
\label{eq_hypmot_sfa}
s = \Cosh \delta \cdot \abs{\vro_{1}-\vro_{2}},
\EE
and the chordlength $ L $ as
\BE
\label{ea_hypmot_Lfa}
\Sinh \frac{L}{2} = 
\Cosh \delta \cdot \Sinh \frac{\abs{\vro_{1}-\vro_{2}}}{2}.
\EE
Eliminating the difference $ \abs{\vro_{1}-\vro_{2}} $ from
these two equations yields the assertion \eqref{eq_hypmot_LS6}.

\smallskip
3) \TI{For the third group \eqref{eq_hypmot_1pg}:}

For an arbitrary initial point $ a \in \H $ one calculates
$ \Dot{u}{a} = (a_{0}-a_{1})^{3} $ and 
$ \Dot{u}{u} = 0 $. By a suitable rotation around $ O $, the
$ u $ can be reduced to $ u = (0,1,0) $.
With the choice $ a = (1,1,0) $ (corresponding to $ O $) 
the orbit is
obtained as noted in \eqref{eq_hypmot_LS7}. It cuts every
line in $ \B $ through $ U $ orthogonally and runs 
through the origin $ O $ (setting $ b = 0 $). So the image
of $ \Gamma $ is part of $ \Hc $. 
In fact, every point on $ \Hc $ is
reached, as follows from the limit relations
$$
\lim_{b\to \pm\infty} 
\frac{\frac{1}{2}b^{2}}{1+\frac{1}{2}b^{2}} = 1,\qquad
\lim_{b\to \pm\infty} 
\frac{-b}{1+\frac{1}{2}b^{2}} = \mp 0.
$$
That each point is reached once is directly clear from the
parametrization.

From \eqref{eq_hypmot_arcl}, the 
arclength on $ \Hc $ between two points 
$ \Gamma(b_{1}) $ and $ \Gamma(b_{2}) $ is calculated to be
\BE
\label{eq_hypmot_sfg}
s = \abs{b_{1}-b_{2}},
\EE
and the chordlength $ L $ as
\BE
\label{eq_hypmot_Lfa}
\Sinh \frac{L}{2} = 
\frac{1}{2}\abs{b_{1}-b_{2}}.
\EE
Again, an elimination yields \eqref{eq_hypmot_LS9}.

The first argument shows that an arbitrary horocycle can be 
transformed by a rotation in such a way that its horizon
point $ U $ becomes $ (1,0) $. For two horocycles with the
same horizon point $ U = (1,0) $, the hyperbolic equivalence 
is proved as follows: The linear maps $ T(\vro) $ from the
second group in \eqref{eq_hypmot_1pg} leave the quadric
shell $ \H $ fixed (in fact the scalar product itself), and 
$ u = (1,1,0) $ is an eigenvector of each $ T(\vro) $,
namely $ T(\vro)u = e^{\vro} u $. For the
image vectors $ y := T(\vro)x $ of $ x $ with 
$ \Dot{u}{x} = p $ one has $ \Dot{T(\vro)u}{y} = p $, thus
$ \Dot{u}{y} = e^{-\vro} p $. So, if $ \Dot{u}{y} = q $ is
the equation of a second horocycle with the same $ u $, this
horocycle is obtained from the first one by a hyperbolic 
translation, represented by $ T(\vro) $, if one chooses 
$ q =  e^{-\vro}p $, i.e. $ \vro = \ln\, p - \ln\, q $. 
\QED

\medskip
The arguments $ \vhi $ in \eqref{eq_hypmot_LS1} are called
\TI{angle values} though it is not necessary here to
understand them as angles in the sense of definition
\eqref{eq_prelim_metr}. More likely, the 
$ \vhi $ should be viewed as parametrizing the universal
covering of $ \Sc(R) $.


\section{\hspace{-1em}.
Closedness conditions}
\label{closed}
\markright{\ref{closed}. Closedness conditions}

A \TI{polygon} in the hyperbolic plane $ \B $ is determined 
by a list of points 
$ Z_{1},\ldots,Z_{n} \in \B $, the \TI{vertices} with 
$ Z_{k} \neq Z_{k+1} $ for $ k = 1,\ldots,n-1 $ and also 
$ Z_{n} \neq Z_{1} $ with $ n \geq 3 $ fixed. Associated to 
these points is the closed \TI{polygon chain} of segments
from $ Z_{k} $ to $ Z_{k+1} $, $ k =1,\ldots,n-1 $ and 
$ Z_{n} $ to $ Z_{1} $. In fact, this polygon chain is the
main object here. The polygon chain remains 
unaltered under cyclic permutations of the vertices. So, two
polygons will be considered the same if they differ 
from each other by a cyclic permutation of the vertices.
Other permutations of the vertices are not allowed. 
The annotation $ Z_{1}\ldots Z_{n} $ denotes the polygon in
this sense.
If nothing else is said, the
indices of vertices will henceforth counted modulo 
$ n $, e.g. $ Z_{n+1} = Z_{1} $. For a vertex $ Z_{k} $,
$ k = 1,\ldots,n $, the two vertices $ Z_{k-1} $ and 
$ Z_{k+1} $ are called \TI{adjacent} to $ Z_{k} $.
For $ k = 1,\ldots,n $, the
lines $ Z_{k} \vee Z_{k+1} $, resp.
the segments $ [Z_{k}Z_{k+1}] $ are called the
\TI{edgelines}, resp. the \TI{edges} of the polygon
(possibly directed if necessary).
Polygons of this type are
fairly general; they allow self-intersections of dimension 
$ 0 $ and even $ 1 $. In order to specify the number $ n $, 
sometimes the term \TI{$ n $-gon} will be used. The
\TI{sidelengths} and the \TI{perimeter} are defined by 
$$
L_{k} := d(Z_{k},Z_{k+1}), \qquad
L := \sum_{k=1}^{n} L_{k}.
$$
The list $ L_{1},\ldots,L_{n} $, corresponding to the list
of vertices $ Z_{1},\ldots,Z_{n} $, is sometimes called a
\TI{length spectrum}. 
Any points $ Z_{1},\ldots,Z_{n} \in \B $ are named 
\TI{collinear}, resp. \TI{cocyclic} if they are situated on a
line resp. circle (in the general sense).  Representing 
vectors $ z_{1},\ldots,z_{k} $ have to
obey the general conventions from Sect. \ref{prelim}. So the
lower indices of the $ z_{k} $ don't count coordinates. 
(If necessary, the coordinates of $ z_{k} $ must be denoted
by $ z_{k0}, z_{k1}, z_{k2} $.)

Convexity in the circle model of Cayley/Klein is not so much
different from vector space convexity because the segments
are the same in both geometries. For polygons, the only
convexity notion to be used here is the following: 
A polygon $ Z_{1}\ldots Z_{n} $ shall be called 
\TI{oriented-convex} if, for any $ k \in \{1,\ldots,n\} $
all vertices $ Z_{j} $ with 
$ j \in \{1,\ldots,n\}\setminus\{k,k+1\} $ lie on the left
hand side of the directed line $ \pfeil{Z_{k} \vee Z_{k+1}} $;
equivalently
\BE
\label{eq_polyrig_pos0}
[Z_{k},Z_{k+1},Z_{j}] > 0 \qquad 
\forall \quad k \in \{1,\ldots,n\}, \quad
j \in \{1,\ldots,n\}\setminus\{k,k+1\}.
\EE
Then, for the same indices, no edgeline $ Z_{k}\vee Z_{k+1} $
contains another vertex $ Z_{j} $, in particular
all vertices $ Z_{1},\ldots,Z_{n} $ are
pairwise distinct. Also, no edgeline cuts another edge in
its interior.

\medskip
The notion of oriented convexity has the advantage that
there is no necessity to specify any `interior' and that 
it can be verified by finitely many inequalities on the
vertices. 

\medskip
\TI{Remark.} 
As to the usual (Euclidean) convexity, one can show that the
polygon chain of an oriented-convex polygon is the boundary
of a unique convex and bounded domain $ C $ whose closure 
$ \overline{C} $ has the vertices as extreme points. So, by
the Krein/Milman theorem, $ \overline{C} $ is the closed
convex hull of the set of vertices. --
Conversely, given the convex hull $ C $ of finitely many
points $ W_{1},\ldots,W_{m} $, the set $ C $ is compact (in
particular $ C = \overline{C} $) and it is the convex hull
of the set $ \{Z_{1},\ldots,Z_{n}\} $ of extreme points of 
$ C $ which is a subset of $ \{W_{1},\ldots,W_{m}\} $. The
extreme points are situated on the boundary of $ C $ which
may be viewed as an oriented simply closed curve. Counted in
the (cyclic) order of this curve, the extreme points form a
polygon as defined in the beginning of this section.

However in this work, all arguments are based on the notion
of oriented convexity thus being independent of the
foregoing relations.

\medskip
Next we consider curved paths which eventually may contain
the vertices of a polygon. Generally, a path in $ \B $ is a
continuous map $ \Gamma: J \to \B $ ($ J $ an interval with
non-void interior). As above, a path may be represented by a
continuous map $ \gamma: J \to \Rd $. 

The path $ \Gamma $ is called a \TI{circum-path} of a polygon
$ Z_{1}\ldots Z_{n} $, if there are arguments 
$ t_{1},\ldots,t_{n} \in J $ with $ Z_{k} = \Gamma(t_{k}) $, 
$ k = 1,\ldots, n $. 

In particular, we shall consider injective paths 
with open interval $ J $ and the following
property: For all $ t_{1} < t_{2} $ in $ J $, the part of 
the image $ \Gamma(J) $ on the right hand side of the
directed line is equal to 
$ \Gamma(\,\left]t_{1},t_{2}\right[\,) $. This part will be 
called the \TI{arc} over the chord 
$ \Gamma(t_{1}),\Gamma(t_{2}) $ and denoted by 
$ B_{\Gamma}(t_{1},t_{2}) $. Such a path shall be called 
a \TI{leftcurve}. Equivalent to this definition is the
requirement 
$ [\gamma(t_{1}),\gamma(t_{2}),\gamma(t)] < 0 $ for all
$ t_{1} < t < t_{2} $ in $ J $ or also
\BE
\label{eq_polyrig_leftc}
[\gamma(t_{1}),\gamma(t),\gamma(t_{2})] > 0 \quad
\text{for all $ t_{1} < t < t_{2} $ in $ J $}.
\EE

\BL
\label{lem_polyrig_equiv}
A leftcurve $ \Gamma: J \to \B $ is the circum-path of an 
oriented-convex $ n $-gon if and only if there is a cyclic
permutation of its vertices $ Z_{1},Z_{2},\ldots,Z_{n} $
such that there are arguments 
$ t_{1} < t_{2} < \cdots < t_{n} $ in $ J $ with 
$ Z_{k} = \Gamma(t_{k}) $, $ k = 1,\ldots,n $.
In this case one has:
\BE
\label{eq_polyrig_equiv1}
B_{\Gamma}(t_{1},t_{n}) \setminus \{Z_{2},\ldots,Z_{n-1}\} =
\Dunl{k=1}{n-1} B_{\Gamma}(t_{k},t_{k+1}).
\EE
\EL

\vspace{-3ex}
The boldface dot signalizes \TI{disjoint} union.

\medskip
\TI{Proof of \ref{lem_polyrig_equiv}.}

\TI{First direction: Assume $ Z_{k} = \Gamma(t_{k}) $ with
$ t_{1} < t_{2} < \cdots < t_{n} $.} 

One has to show 
$ [\gamma(t_{k}),\gamma(t_{k+1}),\gamma(t_{j})] > 0 $ for 
the indices as in \eqref{eq_polyrig_pos0}.

Case $ 1 \leq k \leq n-1 $:
If  $ t_{k} < t_{k+1} < t_{j} $, the assertion is clear from
\eqref{eq_polyrig_leftc}, and likewise if 
$ t_{j} < t_{k} < t_{k+1} $ since \eqref{eq_polyrig_leftc}
implies 
$ [\gamma(t_{j}),\gamma(t_{k}),\gamma(t_{k+1})] > 0 $, hence
$ [\gamma(t_{k}),\gamma(t_{k+1}),\gamma(t_{j})] > 0 $.

Case $ k = n $: 
Then $ t_{k+1} = t_{1} $, and necessarily 
$ t_{1} < t_{j} < t_{n} $, so by \eqref{eq_polyrig_leftc} 
$ [\gamma(t_{1}),\gamma(t_{j}),\gamma(t_{n})] > 0 $, hence
$ [\gamma(t_{n}),\gamma(t_{1}),\gamma(t_{j})] > 0 $.

\smallskip
\TI{Second direction: Assume the polygon to be
oriented-convex in the original arrangement
$ Z_{1},\ldots,Z_{n} $ of vertices.}

In order to line up the vertices monotonically along the
leftcurve let the numbering already be arranged in such a
way that $ Z_{1} $ is the point with the smallest curve
parameter: $ Z_{1} = \Gamma(t_{1} $) and $ t_{1} < t_{j} $
for all $ j = 2,\ldots,n $. Then $ Z_{2} $ has a curve
parameter $ t_{2} $ with $ t_{1} <  t_{2} $. All other 
curve parameters $ t_{j} $, $ j = 3,\ldots,n $ are then above 
$ t_{2} $. For, on account of the oriented convexity of the 
polygon, one has 
$ [\gamma(t_{1}),\gamma(t_{2}),\gamma(t_{j})] > 0 $.
If one had $ t_{1} < t_{j} < t_{2} $ \TI{this} would imply
by the definition of a leftcurve
$ [\gamma(t_{1}),\gamma(t_{j}),\gamma(t_{2})] > 0 $, a
contradiction. By inductively repeating this conclusion 
($ t_{2} $ takes over the role of $ t_{1} $ etc.) it results that
the curve parameters of the vertices can be chosen such that
$ Z_{k} = \Gamma(t_{k}) $ for 
$ k = 1,\ldots,n $ and $ t_{1} < t_{2} < \cdots t_{n} $.

Having achieved this arrangement, Eqn. \eqref{eq_polyrig_equiv1} 
simply follows from the interval decomposition
$$ 
\left]t_{1},t_{n}\right[ =
(\left]t_{1},t_{2}\right[ \dun \{t_{2}\}) \dun
\cdots \dun
(\left]t_{n-2},t_{n-1}\right[ \dun \{t_{n-1}\}) \dun
\left]t_{n-1},t_{n}\right[
$$
by applying the injective map $ \Gamma $.
\QED

\medskip
Lemma \ref{lem_polyrig_equiv} will now be applied to
cocyclic polygons. So $ \Gamma $ will be part of a cone
section, namely part or the whole of a hyperbolic circle.
The circle may be assumed in standard position 
(see \ref{rem_circ_d23}). Any such circle 
is a leftcurve if suitably oriented. However, a distance
circle must be punctured, i.e. one point 
must be removed. This point can still be chosen conveniently.
In any case, the vertices of an oriented-convex and cocyclic
polygon should be numbered such that they correspond to
strictly monotonic increasing parameter values. As 
parametrizations of the circles those of Lemma 
\ref{lem_hypmot_LS} will be used.

\medskip
The following closedness conditions are fundamental for the 
main results. First the cases of distance lines and
horocycles will be discussed since they don't have to be
punctured.

\medskip
\BL[distance lines and horocycles]
\label{lem_polyrig_agbil}
Let $ Z_{1}\ldots Z_{n} $ be an oriented-convex polygon and 
without loss of generality let $ L_{n} $ denote its maximal 
sidelength.

If the vertices lie on the distance line $ \Dc(\delta) $,
resp. on the horocycle $ \Hc $ the following 
\underline{closedness} \underline{conditions} hold true:
\begin{align}
\label{eq_polyrig_agbil1}
\Psi(\delta,L_{n}) &= 
\Psi(\delta,L_{1})+\cdots+\Psi(\delta,L_{n-1})
\intertext{resp.}
\label{eq_polyrig_agbil2}
\Omega(L_{n}) &= 
\Omega(L_{1})+\cdots+\Omega(L_{n-1}),
\end{align}
where the functions $ \Psi $ and $ \Omega $ are defined in 
Lemma \ref{lem_hypmot_LS}. 

In both cases the edge of maximal length is unique. 
\EL

\Bew
Eqn. \eqref{eq_polyrig_equiv1} implies for
the corresponding arclengths:
$ s_{n} = s_{1}+\cdots+s_{n-1} $. Using the relations to
the sidelengths from Lemma
\ref{lem_hypmot_LS} one immediately deduces 
\eqref{eq_polyrig_agbil1} and \eqref{eq_polyrig_agbil2}.
The terms on the right hand side of these equations are all 
positive, so 
$ \Psi(\delta,L_{n}) > \Psi(\delta,L_{k}) $, resp. 
$ \Lambda(\delta,L_{n}) > \Lambda(\delta,L_{k}) $
for $ k = 1,\ldots,n-1 $. By the
strict monotony of $ \Omega $ and $ \Lambda $: 
$ L_{n} > L_{k} $.
\QED

\medskip
\BC
\label{cor_polyrig_noreg}
Regular oriented-convex polygons whose vertices are on a
distance line or a horocycle do not exist.
\QED
\EC

Now, in case of a distance circle as circum-path for an
oriented-convex polygon $ Z_{1}\ldots Z_{n} $ there is an
additional peculiarity. Namely there are principally two different
positions of the vertices relative to the center $ C $.
(Here, we may assume $ C = O $ with the point vector 
$ c = (1,0,0) $.)

\smallskip
\TI{Case I:}
There exists a directed edgeline, leaving the circle center 
$ C $ on its right hand side, i.e.
\begin{align*}
\label{eq_polyrig_nische}
\tag{N}
\Exists k \in \{1,\ldots,n-1\}: \quad
[z_{k},z_{k+1},c] < 0
\qquad &\text{(\TI{niche position})}.
\intertext{\TI{Case II:} \TI{Case I is not satisfied, i.e.}}
\label{eq_polyrig_voll}
\tag{F}
\Forall k \in \{1,\ldots,n-1\}: \quad 
[z_{k},z_{k+1},c] \geq 0
\qquad &\text{(\TI{full position})}.
\end{align*}
The names spring from the fact that in case I there is a
closed half circle disk clear of vertices while in case II
every closed half circle disk contains vertices of the
polygon. See the following pictures and Lemmas
\ref{lem_polyrig_dnbil} and \ref{lem_polyrig_dvbil}.

\begin{center}
\begin{minipage}{12.5cm}
\bild{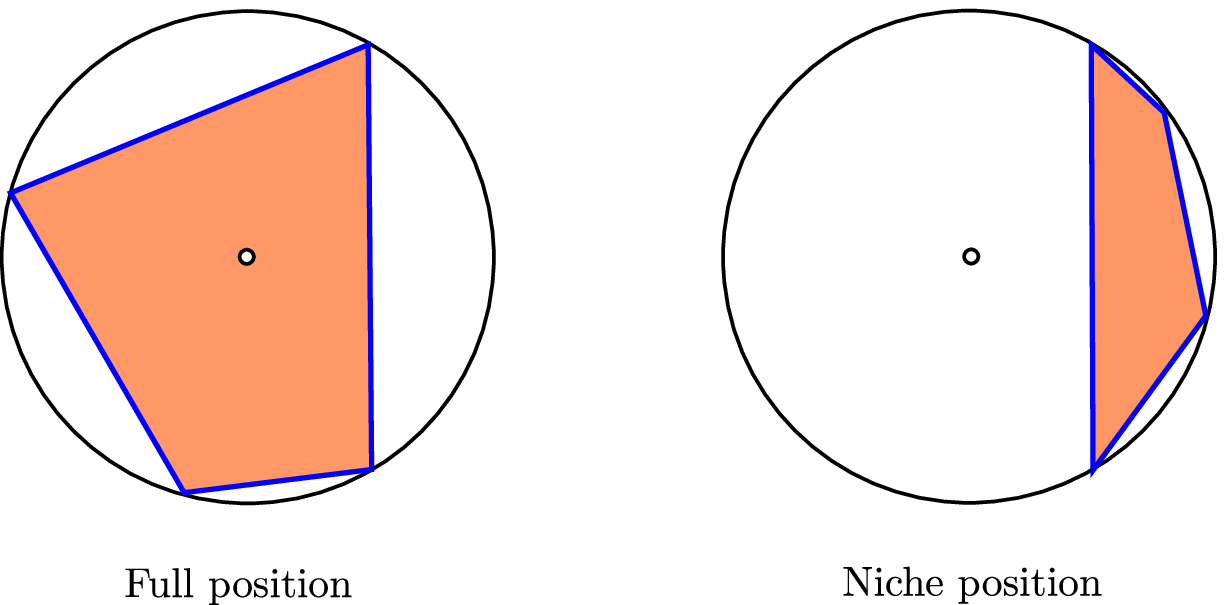}
\end{minipage}
\end{center}

\BL
\label{lem_polyrig_numm}
Let an oriented-convex $ n $-gon have the distance circle 
$ \Sc(R) $ as an circum-path. Then it is possible to
arrange the numbering of vertices, the puncturing of 
$ \Sc(R) $, and a rotation around the center $ O $ such that
an arbitrarily chosen edgeline becomes 
$ \pfeil{Z_{n} \vee Z_{1}} $ and such that 
$ Z_{k} = \Gamma(\vhi_{k}) $ with
$ 0 < \vhi_{1} < \vhi_{2} < \cdots \vhi_{n} < 2 \pi $.
\EL

A puncturing of $ \Sc(R) $ corresponds to an angle
determination of the circle points according to 
\eqref{eq_hypmot_LS1} in that the removed point obtains an
angle value $ \nu_{0} $ and all other circle points receive 
angle values $ \vhi $ in the open interval 
$ \left]\nu_{0},\nu_{0}+2\pi\right[ $. Of course, each
distance circle is to be oriented such that, after
puncturing, a leftcurve arises. 

\medskip
\TI{Proof of \ref{lem_polyrig_numm}.}
First, as a point to remove, one can take any $ P \in S(R) $
different from a vertex. Then one may rotate the whole
figure such that $ P $ lands on the positive part of the
groundline. Further, by Lemma \ref{lem_polyrig_equiv}, the
vertices can be numbered such that, in the
parametrization \eqref{eq_hypmot_LS1}, there holds
$ z_{k} = \gamma(\vhi_{k}), k = 1,\ldots,n $ with
$ 0 < \vhi_{1} < \vhi_{2} < \cdots \vhi_{n} < 2\pi $.
Now let $ \pfeil{Z_{j}\vee Z_{j+1}} $ be any edgeline with 
$ j = 1,\ldots n-1 $ (for $ j = n $ the assertion is
already valid). Then, for the vertices starting with 
$ Z_{j+1} $, there persists the following list with possible new
angle values in the second line and possible new
names for the vertices in the third line:
$$
\begin{array}{ccccccccccccccc}
Z_{j+1}    &  & Z_{j+2}    &  & \dots &  & Z_{n}    & 
 & Z_{1} &  & Z_{2} &  & \dots   &  & Z_{j} 
\\[1ex]
\vhi_{j+1} & < & \vhi_{j+2} & < & \dots & < & \vhi_{n} &
< & \vhi_{1}+2\pi & < & \vhi_{2}+2\pi & < & \cdots  & < & \vhi_{j}+2\pi 
\\[1ex]
Z_{1}'    &  & Z_{2}'      &  & \dots &  & Z_{n-j}'  & 
 & Z_{n-j+1}' &  & Z_{n-j+2}'   &  & \dots   &  & \;Z_{n}'.
\end{array}
$$
One already has
$ \pfeil{Z_{n}'\vee Z_{1}'} = \pfeil{Z_{j}\vee Z_{j+1}} $.
In the second line, the new angle values are strictly
monotonic increasing, in particular at 
$ \vhi_{n} < \vhi_{1}+2\pi $, since this is equivalent to 
$ \vhi_{n} -\vhi_{1} < 2\pi $. For the difference of the 
`ends' one has $ (\vhi_{j}+2\pi)-\vhi_{j+1} < 2\pi $. So 
there is a new angle determination for which the vertices in
the sequence of the third line are represented by strictly 
monotonic increasing angle values. By another rotation
this angle determination can be reduced to the interval 
from $ 0 $ to $ 2\pi $.
\QED

\medskip
\BL[distance circles -- niche position]
\label{lem_polyrig_dnbil}
Let an oriented-convex $ n $-gon have the distance circle 
$ \Sc(R) $ as an circum-path and assume it in niche-situation
there. 
Without loss of generality, let $ \pfeil{Z_{n}\vee Z_{1}} $
be an edgeline leaving the center $ O $ on its right hand
side. Then
\begin{DES1}{(iii)}
\item[(i)]
The open halfplane right of $ \pfeil{Z_{n}\vee Z_{1}} $
doesn't contain vertices. A fortiori, there are closed
halfcircle disks clear of vertices.
\item[(ii)]
All the edgelines $ \pfeil{Z_{k}\vee Z_{k+1}} $,
$ k = 1,\ldots,n-1 $, leave the center $ O $ on their left
hand side. So, there is only one edgeline,
$ \pfeil{Z_{n}\vee Z_{1}} $, having $ O $ on its right hand 
side.
\item[(iii)]
For the sidelengths, there holds the 
\underline{closedness condition}
\BE
\label{eq_polyrig_dnbil1}
\Phi(R,L_{n}) = \Phi(R,L_{1})+\cdots+\Phi(R,L_{n-1}),
\EE
where the function $ \Phi $ is defined in 
Lemma \ref{lem_hypmot_LS}. 
\item[(iv)]
The edgeline $ \pfeil{Z_{n}\vee Z_{1}} $ is also
characterized by belonging to the biggest side of the
polygon:
$$
d(Z_{1},Z_{n}) > d(Z_{k},Z_{k+1}), \quad
k = 1,\ldots,n-1. 
$$
\end{DES1}
\EL

\Bew
As a distinguished edgeline $ \pfeil{Z_{n}\vee Z_{1}} $ 
in the sense of Lemma \ref{lem_polyrig_numm} one can choose 
one which has the center $ C = O $ on its right hand side, so
$ [z_{n},z_{1},c] < 0 $. By
\eqref{eq_hypmot_LS1} this means
$$
0 < [z_{1},z_{n},c] =
\BV
\;\Cosh R & \Cosh R & 1\; \\
\;\Sinh R \Cos \vhi_{1} & \Sinh R \Cos \vhi_{n} & 0\; \\
\;\Sinh R \Sin \vhi_{1} & \Sinh R \Sin \vhi_{n} & 0\; 
\EV
= \sinh^{2}R\cdot \sin(\vhi_{n}-\vhi_{1}),
$$
hence $ 0 < \vhi_{n}-\vhi_{1} < \pi $.

\smallskip
\TI{For (i):} 
All other vertices lie left of $ \pfeil{Z_{n}\vee Z_{1}} $. 
So right of $ \pfeil{Z_{n}\vee Z_{1}} $ there are no points 
of the polygon chain. Then there exist directed lines
through $ C $ whose closed right hand sides are clear of such
points, e.g. the line through $ C $ parallel in the Euclidean 
sense and equally directed to $ \pfeil{Z_{n}\vee Z_{1}} $.

\smallskip
\TI{For (ii):}
It must be shown 
$ [\gamma(\vhi_{k}),\gamma(\vhi_{k+1}),c] > 0 $ for 
$ k = 1,\ldots,n-1 $. Now, as above
$
[\gamma(\vhi_{k}),\gamma(\vhi_{k+1}),m] 
= \sinh^{2}R\cdot \sin(\vhi_{k+1}-\vhi_{k})
$,
and, since 
$ 0 < \vhi_{k+1}-\vhi_{k} < \vhi_{n}-\vhi_{1} < \pi $, the
assertion follows.

\smallskip
\TI{For (iii):} 
Again, this is deduced from \eqref{eq_polyrig_equiv1},
the definition of $ B_{\Gamma} $, and Eqn. 
\eqref{eq_hypmot_LS3} (first part).

\smallskip
\TI{For (iv):} 
From (iii) follows, because all terms on the right hand side
are positive, that 
$ \Phi(R,L_{n}) > \Phi(R,L_{k}) $ for $ k = 1,\ldots,n-1 $.
Then the strict monotony of the function $ \Phi $ in its
second argument implies $ L_{n} > L_{k} $.
\QED

\medskip
\BL[distance circles -- full position]
\label{lem_polyrig_dvbil}
Let an oriented-convex $ n $-gon have the distance circle 
$ \Sc(R) $ as an circum-path and assume it in full position
there. Then:
\begin{DES1}{(ii)}
\item[(i)]
Each half circle disk contains vertices of the polygon.
\item[(ii)]
For the sidelengths, there holds the 
\underline{closedness condition}
\BE
\label{eq_polyrig_dvbil2}
2\pi\Sinh R = \Phi(R,L_{1})+\cdots+\Phi(R,L_{n}).
\EE
\end{DES1}
\EL

\Bew

\TI{For (i):}
Assume that a half circle disk would contain no vertices.
Then, by an eventual rotation around $ O $ one could achieve
$ z_{k} = \gamma(\vhi_{k}) $ with
$ 0 < \vhi_{1} < \cdots \vhi_{n} < \pi $. By Lemma 
\ref{lem_polyrig_dnbil} (i) this would imply 
$ [z_{1},z_{n},c] > 0 $, so $ [z_{n},z_{1},c] < 0 $. The
edgeline $ \pfeil{Z_{n}\vee Z_{1}} $ then had $ O $ on its
right hand side.

\smallskip
\TI{For (ii):}
Let everything be arranged like in Lemma
\ref{lem_polyrig_numm}, where $ \pfeil{Z_{n}\vee Z_{1}} $
may be any edgeline, e.g. one with maximal sidelength.
Again, by  \eqref{eq_polyrig_equiv1}:
\BE
\label{eq_polyrig_dvbil1}
B_{\Gamma}(\vhi_{1},\vhi_{n}) \setminus \{Z_{2},\ldots,Z_{n-1}\} =
B_{\Gamma}(\vhi_{1},\vhi_{2}) 
\dun \cdots \dun 
B_{\Gamma}(\vhi_{n-1},\vhi_{n}).
\EE
All the arcs on the right hand side of 
\eqref{eq_polyrig_dvbil1} have the center $ O $ on their
left hand sides since $ O $ is always left of 
$ \pfeil{Z_{k}\vee Z_{k+1}} $ and the arc 
$ B_{\Gamma}(\vhi_{1},\vhi_{n}) $ is always right of 
$ \pfeil{Z_{k}\vee Z_{k+1}} $. So, for the relation between 
arc- and chordlength, the first part of
\eqref{eq_hypmot_LS3} applies: $ s_{k} = \Phi(R,L_{k}) $, $
k = 1,\ldots,n-1 $.

For the arc on the left hand side of
\eqref{eq_polyrig_dvbil1}, the situation is converted: This 
arc has $ O $ on the right hand side since $ O $ lies left
of $ \pfeil{Z_{n}\vee Z_{1}} $. So here, the relation 
between arc- and chordlength is regulated by the second part
of \eqref{eq_hypmot_LS3}: $ s_{n} = 2\pi\Sinh R-\Phi(R,L_{n}) $. 

From \eqref{eq_polyrig_dvbil1} follows 
$ s_{n} = s_{1}+\cdots+s_{n-1} $. With the above values of the
sidelengths, Eqn. \eqref{eq_polyrig_dvbil2} is thus verified.
\QED


\section{\hspace{-1em}.
Rigidity}
\label{rigi}
\markright{\ref{rigi}. Rigidity}

A first step to the rigidity of oriented-convex cocyclic
polygons is the uniqueness of the circle type to be
formulated in Theorem \ref{thm_polyrig_utyp}. The
following preparations are needed:

\medskip
\BL
\label{lem_polyrig_Sinsin}
For positive real numbers $ x_{1},\ldots,x_{m} $ with
$ m \geq 2 $ there holds:
\BE
\label{eq_polyrig_Sinsin1}
\sinh(x_{1}+\cdots+x_{m}) 
> \Sinh x_{1}+\cdots+\Sinh x_{m}
\EE
and, if $ x_{1}+\cdots+x_{m} < \pi $:
\BE
\label{eq_polyrig_Sinsin2}
\sin(x_{1}+\cdots+x_{m}) 
< \Sin x_{1}+\cdots+\Sin x_{m}.
\EE
\EL

\Bew
In both cases one can proceed inductively:

\TI{For \eqref{eq_polyrig_Sinsin1}:}
If $ m = 2 $ then
$$
\sinh(x_{1}+x_{2}) = 
\Sinh x_{1} \Cosh x_{2}+\Sinh x_{2} \Cosh x_{1} > 
\Sinh x_{1}+\Sinh x_{2}.
$$
The induction step from $ m $ to $ m+1 $ is deduced from
$$
\sinh(x_{1}+\cdots+x_{m+1}) >
\sinh(x_{1}+\cdots+x_{m})+\Sinh x_{m+1} >
\Sinh x_{1}+\cdots+\Sinh x_{m}+\Sinh x_{m+1}.
$$

\TI{For \eqref{eq_polyrig_Sinsin2}:}
If $ m = 2 $ then
$$
\sin(x_{1}+x_{2}) = 
\Sin x_{1} \Cos x_{2}+\Sin x_{2} \Cos x_{1} < 
\Sin x_{1}+\Sin x_{2}.
$$
The induction step from $ m $ to $ m+1 $ is deduced from
$$
\sin(x_{1}+\cdots+x_{m+1}) <
\sin(x_{1}+\cdots+x_{m})+\Sin x_{m+1} <
\Sin x_{1}+\cdots+\Sin x_{m}+\Sinh x_{m+1}
$$
where of course the additional assumption enters.
\QED

\medskip
\BR
\label{rem_polyrig_umk}
Transcribing Eqns. \eqref{eq_polyrig_Sinsin1} and 
\eqref{eq_polyrig_Sinsin2} to the inverse functions
yields for positive real numbers $ \xi_{1},\dots,\xi_{m} $
with $ m \geq 2 $:
\begin{align}
\label{rem_polyrig_umk1}
\arsinh(\xi_{1}+\cdots+\xi_{m}) &< 
\Arsinh \xi_{1}+\cdots+\Arsinh \xi_{m} \\
\label{rem_polyrig_umk2}
\arcsin(\xi_{1}+\cdots+\xi_{m}) &>
\Arcsin \xi_{1}+\cdots+\Arcsin \xi_{m}, \quad
\xi_{1}+\cdots+\xi_{m} \leq \frac{\pi}{2}.
\end{align}
\ER

\smallskip
\BT
\label{thm_polyrig_utyp}
Solely by the sidelengths of a cocyclic oriented-convex
polygon it is determined on which type of circle the
vertices are situated.
\ET

\Bew
Decisive for this are the closedness conditions from 
Lemmas \ref{lem_polyrig_agbil}, 
\ref{lem_polyrig_dnbil}, and
\ref{lem_polyrig_dvbil}. It will be shown that
any two of these conditions exclude each other.

\smallskip
\TI{For the pair distance line/horocycle:}

Assume that at the same time the following would be true:
\begin{align}
\label{eq_polyrig_utyp1}
\Psi(\delta,L_{n}) &= 
\Psi(\delta,L_{1})+\cdots+\Psi(\delta,L_{n-1}) \\[0.5ex]
\label{eq_polyrig_utyp2}
\Omega(L_{n}) &= 
\Omega(L_{1})+\cdots+\Omega(L_{n-1}).
\end{align}
Then, eliminating the expression $ \Sinh\frac{t}{2} $ from Eqns. 
\eqref{eq_hypmot_LS5}, \eqref{eq_hypmot_LS8} yields
$$
\frac{\Omega(t)}{\delta_{1}} = 
\Sinh\frac{\Psi(\delta,t)}{\delta_{1}},
\qquad \delta_{1} := 2\Cosh\delta.
$$
From \eqref{eq_polyrig_utyp1} and \eqref{eq_polyrig_utyp2}
follows by setting $ t = L_{n} $, resp. $ t = L_{j} $:
$$
\Sinh\frac{\Psi(\delta,L_{1})}{\delta_{1}}+\cdots+
\Sinh\frac{\Psi(\delta,L_{n-1})}{\delta_{1}}
=
\Sinh\frac{\Psi(\delta,L_{1})+\cdots+
\Psi(\delta,L_{n-1})}{\delta_{1}}.
$$
By \eqref{eq_polyrig_Sinsin1}, such an equation can never hold.

\smallskip
\TI{For the pair distance circle (niche position)/horocycle:}

Assume that at the same time the following would be true:
\begin{align}
\label{eq_polyrig_utyp3}
\Phi(R,L_{n}) &= 
\Phi(R,L_{1})+\cdots+\Phi(R,L_{n-1}) \\[0.5ex]
\label{eq_polyrig_utyp4}
\Omega(L_{n}) &= 
\Omega(L_{1})+\cdots+\Omega(L_{n-1}).
\end{align}
Then, eliminating the expression $ \Sinh\frac{t}{2} $ from Eqns. 
\eqref{eq_hypmot_LS2}, \eqref{eq_hypmot_LS8} yields
$$
\frac{\Omega(t)}{R_{1}} = 
\Sin\frac{\Phi(R,t)}{R_{1}},
\qquad
R_{1} := 2\Sinh R.
$$
From \eqref{eq_polyrig_utyp3} and \eqref{eq_polyrig_utyp4} 
follows by setting $ t = L_{n} $, resp. $ t = L_{j} $:
$$
\Sin\frac{\Phi(R,L_{1})}{R_{1}}+\cdots+
\Sin\frac{\Phi(R,L_{n-1})}{R_{1}}
=
\Sin\frac{\Phi(R,L_{1})+\cdots+
\Phi(R,L_{n-1})}{R_{1}}.
$$
By \eqref{eq_polyrig_Sinsin2}, such an equation can never hold.

\smallskip
\TI{For the pair distance circle (niche position)/distance
line:}

Assume that at the same time the following would be true:
\begin{align}
\label{eq_polyrig_utyp5}
\Phi(R,L_{n}) &= 
\Phi(R,L_{1})+\cdots+\Phi(R,L_{n-1}) \\[0.5ex]
\label{eq_polyrig_utyp6}
\Psi(\delta,L_{n}) &= 
\Psi(\delta,L_{1})+\cdots+\Psi(\delta,L_{n-1}).
\end{align}
Then, eliminating the expression $ \Sinh\frac{t}{2} $ from Eqns. 
\eqref{eq_hypmot_LS2}, \eqref{eq_hypmot_LS5} yields, using
the above shortcuts $ \delta_{1} $ and $ R_{1} $:
$$
\Sinh \frac{\Psi(\delta,t)}{\delta_{1}} =
\frac{R_{1}}{\delta_{1}} \Sin \frac{\Phi(R,t)}{R_{1}}.
$$
From \eqref{eq_polyrig_utyp5} and \eqref{eq_polyrig_utyp6}  
follows by setting $ t = L_{n} $, resp. $ t = L_{j} $:
$$
\Sinh
\frac{\Psi(\delta,L_{1})+\cdots+\Psi(\delta,L_{n-1})}{\delta_{1}}
=
\frac{R_{1}}{\delta_{1}}
\Sin\frac{\Phi(R,L_{1})+\cdots+\Phi(R,L_{n-1})}{R_{1}},
$$
so with \eqref{eq_polyrig_Sinsin1}
$$
\Sinh
\frac{\Psi(\delta,L_{1})}{\delta_{1}}+\cdots+
\frac{\Psi(\delta,L_{n-1})}{\delta_{1}}
<
\frac{R_{1}}{\delta_{1}}
\Sin\frac{\Phi(R,L_{1})+\cdots+\Phi(R,L_{n-1})}{R_{1}},
$$
or
$$
\frac{R_{1}}{\delta_{1}}\Sin\frac{\Phi(R,L_{1}))}{R_{1}}
+\cdots+
\frac{R_{1}}{\delta_{1}}\Sin\frac{\Phi(R,L_{n-1}))}{R_{1}} <
\frac{R_{1}}{\delta_{1}}
\Sin\frac{\Phi(R,L_{1})+\cdots+\Phi(R,L_{n-1})}{R_{1}}.
$$
By \eqref{eq_polyrig_Sinsin2}, such an equation can never hold.

\smallskip
\TI{For the remaining pairs:} 

Now, the last two pairs are lacking, namely replacing the
niche positions with the distance circles by the full
positions. However, the difference in arguing are small.
Only the condition  \eqref{eq_polyrig_utyp3}, resp. 
\eqref{eq_polyrig_utyp5}: 
$ \Phi(R,L_{n}) = \Phi(R,L_{1})+\cdots+\Phi(R,L_{n-1}) $,
has to be replaced by 
$ \pi R_{1} = \Phi(R,L_{1})+\cdots+\Phi(R,L_{n}) $
(see \eqref{eq_polyrig_dnbil1}), i.e. by
$$
\frac{\Phi(R,L_{n})}{R_{1}} =
\pi-
\left(\frac{\Phi(R,L_{1})}{R_{1}}+\cdots+
\frac{\Phi(R,L_{n-1})}{R_{1}}\right).
$$
Since the  sinus value of the right hand side is the same as 
the sinus value of the sum in the last big parenthesis
there result the same inequalities and arguments as with the
last two cases.
\QED

\medskip
The simplest rigidity situation prevails if the vertices lie
on a horocycle:

\medskip
\BT
\label{thm_polyrig_gkongr}
Any two oriented-convex $ n $-gons with vertices on a
horocycle and with the same length spectrum 
$ L_{1},\ldots,L_{n} $ are hyperbolically equivalent.
\ET

\Bew
Again it suffices to consider the special horocycle 
$ \Hc $ from Lemma \ref{lem_hypmot_LS} (iii). 
Change the parameter $ b $ constant proportionally and
equally directed to the arclength parameter $ s $, and let the
corresponding one parameter group be parametrized by 
$ G(s) $ and the horocycle $ \Hc $ by $ \Gamma(s) := G(s)A $.
Between the arclength and the chordlength, one has the
bijective relation expressed by
\eqref{eq_hypmot_LS8}. Permute the length 
numbers such that $ L_{n} $ is the (unique) maximum.

Now, if $ Z_{1} = \Gamma(s_{0}) $ is the first vertex of
such a polygon, then the succeeding vertices necessarily are
given by 
\BE
\label{eq_polyrig_gkongr1}
Z_{k} =
\Gamma(s_{0}+\Omega(L_{1})+\cdots+\Omega(L_{k-1})),
\quad k = 2,\ldots,n;
\EE
see Lemma \ref{lem_polyrig_equiv}.
In other words: The succeeding vertices arise from $ Z_{1} $ 
by laying
off the chordlengths (more accurately: the corresponding
arclengths) along $ \Hc $ in the direction of its
orientation. 

Since 
$ \Gamma(s_{0}+\Omega(L_{1})+\cdots+\Omega(L_{k-1})) = 
G(\Omega(L_{1})+\cdots+\Omega(L_{k-1}))Z_{1} $, these
points only depend on $ Z_{1} $ and the length values 
$ L_{1},\ldots,L_{n-1} $.
Now, if $ Z_{1}' = \Gamma(s_{0}') $ is the first vertex of
another such polygon then Eqn. \eqref{eq_polyrig_gkongr1}
holds correspondingly for its vertices 
$ Z_{2}',\ldots, Z_{n}' $, and from 
$ Z_{1}' = G(s_{0}'-s_{0})Z_{1} $ follows 
$ Z_{k}' = G(s_{0}'-s_{0})Z_{k} $, $ k = 1,\ldots,n $, hence
the equivalence of the two polygons.
\QED

\medskip
This procedure of laying off is also possible for the
remaining cases \TI{insofar} the chordlengths determine the 
corresponding arclength uniquely. But in contrast to the
horocycles the other circle types have invariants, and it
must be clarified in addition that these invariants are
again determined uniquely by the length spectrum.
The following lemmas serve this purpose.

\medskip
\BL
\label{lem_polyrig_ch}
For $ q > 1 $, the function $ h: \Rplusnull \to \R $,
$$
h(x) := \frac{\arsinh(q\Sinh x)}{x},\qquad h(0) := q,
$$
is strictly monotonic decreasing.
In particular
$ h(x) < h(0) = q $ for $ x > 0 $.
\EL

\Bew
The idea is that the values of $ h $ are certain slopes of a
concave function. 

At any rate, the function $ h $ has in $ 0 $ a removable
singularity. With the given value, it is of class 
$ C^{\infty} $ on the whole of $ \Rplusnull $. The numerator
function $ H: \Rplusnull \to \R $, 
$ H(x) := \arsinh(q\Sinh x) $ is strictly concave: Its 
second derivative works out to be
$$
H''(x) =
-q(q^{2}-1)\;\frac{\Sinh x}{(1+q^{2}\sinh^{2}x)^{3/2}},
$$
thus $ H''(x) < 0 $ in $ \Rplus $. For the monotony of the
difference quotients $ \delta H $ this implies 
$ (\delta H)(0,x_{1}) > (\delta H)(0,x_{2}) $ for 
$ 0 < x_{1} < x_{2} $, i.e.
$$
\frac{\arsinh(q\Sinh x_{1})}{x_{1}} >
\frac{\arsinh(q\Sin x_{2})}{x_{2}} \qquad
\text{for $ 0 < x_{1} < x_{2} $},
$$
so $ h(x_{1}) > h(x_{2}) $. From this follows the strictly
monotonic decreasing of $ h $ on the whole of $ \Rplusnull $.
\QED

\medskip
\BC
\label{cor_polyrig_mon}
For $ 0 < \delta' < \delta $ and $ t > 0 $ holds:
\BE
\label{eq_polyrig_quot}
\frac{\Psi(\delta,t)}{\Psi(\delta',t)} > 1
\EE
and also the implication:
\BE
\label{eq_polyrig_quotrel}
0 < t_{1} < t_{2} \Impl 
\frac{\Psi(\delta,t_{1})}{\Psi(\delta',t_{1})} < 
\frac{\Psi(\delta,t_{2})}{\Psi(\delta',t_{2})}.
\EE
\EC

Intuitively, the inequality \eqref{eq_polyrig_quot} says: 
If the same chordlength $ L $ is laid off along the 
distance line $ \Dc(\delta) $ and also along the distance line 
$ \Dc(\delta)' $ with smaller parameter $ \delta' < \delta $ 
then the corresponding arc also becomes smaller: 
$ \Psi(\delta',L) < \Psi(\delta,L) $. And the implication 
\eqref{eq_polyrig_quotrel} means that the arcs over smaller 
chords are less shortened than the arcs over bigger chords
if one passes over from $ \Dc(\delta) $ to $ \Dc(\delta') $.

\medskip
\TI{Proof of \ref{cor_polyrig_mon}.}
Set
$$
\delta_{1}:= 2\Cosh \delta, \qquad 
\delta_{1}' := 2\Cosh \delta', \qquad 
q := \frac{\delta_{1}}{\delta_{1}'} > 1, \qquad 
\tau_{1} := \frac{\ds 2\Sinh \frac{t_{1}}{2}}{\delta_{1}}, \qquad 
\tau_{2} := \frac{\ds 2\Sinh \frac{t_{2}}{2}}{\delta_{1}}.
$$
Then, for $ t = t_{2} $, the assertion \eqref{eq_polyrig_quot} 
reads as 
$$
q\;\frac{\Arsinh \tau_{2}}{\arsinh(q\tau_{2})} > 1,
$$
or, with $ x_{2} := \Arsinh \tau_{2} $, as
$$
\frac{\arsinh(q\Sinh x_{2})}{qx_{2}} < 1.
$$
But this is correct by Lemma \ref{lem_polyrig_ch}, last part.

Analogously, the assertion \eqref{eq_polyrig_quotrel} reads 
as
$$
\frac{\arsinh(q\tau_{1})}{\Arsinh \tau_{1}} <
\frac{\arsinh(q\tau_{2})}{\Arsinh \tau_{2}},
$$
or, with  $ x_{1} := \Arsinh \tau_{1} $, as
$$
\frac{\arsinh(q\Sinh x_{1})}{x_{1}} > 
\frac{\arsinh(q\Sinh x_{2})}{x_{2}}.
$$
Since $ 0 < x_{1} < x_{2} $ this is true by Lemma 
\ref{lem_polyrig_ch}.
\QED

\medskip
\BT
\label{thm_polyrig_akongr}
Any two oriented-convex $ n $-gons with vertices on 
distance lines and with the same length spectrum 
$ L_{1},\ldots,L_{n} $ are hyperbolically equivalent.
\ET

\Bew
The main point is to identify the distance invariants 
$ \delta, \delta' $ of two possible circum-distance-lines. 

Again, let the numbering of the sidelengths be arranged such
that $ L_{n} $ is their (unique) maximum, so 
$ L_{n} > L_{k} $ for $ k = 1,\ldots,n-1 $.

Assuming $ \delta' < \delta $, one arrives at a
contradiction in the following way. The two closedness
conditions sound
\begin{align}
\label{eq_polyrig_akongr1}
\Psi(\delta,L_{n}) &= 
\Psi(\delta,L_{1})+\cdots+\Psi(\delta,L_{n-1}) \\[0.5ex]
\label{eq_polyrig_akongr2}
\Psi(\delta',L_{n}) &= 
\Psi(\delta',L_{1})+\cdots+\Psi(\delta',L_{n-1}).
\end{align}
By Corollary \ref{cor_polyrig_mon} 
\BE
\label{eq_polyrig_akongr3}
\frac{\Psi(\delta,L_{k})}{\Psi(\delta',L_{k})} < 
\frac{\Psi(\delta,L_{n})}{\Psi(\delta',L_{n})} =: p,
\qquad k = 1,\ldots,n-1.
\EE
This implies
$$
\frac{\Psi(\delta,L_{n})}{\Psi(\delta',L_{n})} =
\frac{\Psi(\delta,L_{1})+\cdots+\Psi(\delta,L_{n-1})}%
{\Psi(\delta',L_{1})+\cdots+\Psi(\delta',L_{n-1})} <
\frac{p\Psi(\delta',L_{1})+\cdots+p\Psi(\delta',L_{n-1})}%
{\Psi(\delta',L_{1})+\cdots+\Psi(\delta',L_{n-1})} = p,
$$
contradicting \eqref{eq_polyrig_akongr3}.

Now, it can be assumed that the two polygons have the same 
circum-distance-line, e.g. $ \Dc(\delta) $ from 
\ref{rem_circ_d23}. Then the above procedure of 
laying off the chordlengths along $ \Dc(\delta) $ again leads
to the conclusion that the vertices of the first polygon 
are transformed to those of the second polygon by a
common element of the corresponding one parameter group.
\QED

\medskip
For distance circles as circum-paths one may proceed in an
analogous way, where however the two different positions 
(niche, resp. full) must be regarded. Again, the main point 
is to distil out the radius invariant solely from the length
spectrum. For this there are statements completely analogous 
to Lemma \ref{lem_polyrig_ch} and Corollary 
\ref{cor_polyrig_mon}, with analogous proofs. Therefore, the
arguments will be shortened somewhat.

\medskip
\BL
\label{lem_polyrig_dch}
For $ q > 1 $, the function
$$
\begin{aligned}
h&: \left[0,\Arcsin\frac{1}{q}\right] \To \R \\
h(x) &:= \frac{\arcsin(q\Sin x)}{x},\qquad h(0) := q
\end{aligned}
$$
is strictly monotonic increasing.
In particular
$ q = h(0) < h(x) $ for $ \ds 0 < x \leq \Arcsin\frac{1}{q} $.
\EL

\Bew
As above, this runs with the aid of the numerator function 
$ H(x) := \arcsin(q\Sin x) $ because $ h(x) $ represents the
slope of $ H $ between $ 0 $ and $ x $. This time, $ H $ is 
convex with the effect that the slopes are strictly
monotonic increasing.
\QED

\medskip
\BC
\label{cor_polyrig_dmon}

For $ 0 < R' < R $ and $ 0 < t \leq 2R' $ holds:
\BE
\label{eq_polyrig_dquot}
\frac{\Phi(R,t)}{\Phi(R',t)} < 1
\EE
and also the implication: 
\BE
\label{eq_polyrig_dquotrel}
0 < t_{1} < t_{2} \leq 2R' \Impl 
\frac{\Phi(R,t_{1})}{\Phi(R',t_{1})} > 
\frac{\Phi(R,t_{2})}{\Phi(R',t_{2})}.
\EE
\EC

\Bew
One can perform the same steps as for Corollary 
\ref{cor_polyrig_mon}. This time, setting
$$
R_{1}:= 2\Sinh R, \qquad 
R_{1}':= 2\Sinh R', \qquad 
q := \frac{R_{1}}{R_{1}'} > 1, \qquad 
\tau_{1} := \frac{\ds 2\Sinh \frac{t_{1}}{2}}{R_{1}}, \qquad 
\tau_{2} := \frac{\ds 2\Sinh \frac{t_{2}}{2}}{R_{1}},
$$
the assertion  \eqref{eq_polyrig_dquot} reads for 
$ t = t_{2} $ as
$$
q\;\frac{\Arcsin \tau_{2}}{\arcsin(q\tau_{2})} < 1 
$$
or, with  $ x_{2} := \Arcsin \tau_{2} $, as
$$
\frac{\arcsin(q\sin x_{2})}{qx_{2}} > 1.
$$
This is true by Lemma \ref{lem_polyrig_dch}, last part.

Similarly, the assertion  \eqref{eq_polyrig_dquotrel} reads 
as
$$
\frac{\arcsin(q\tau_{1})}{\Arcsin \tau_{1}} <
\frac{\arcsin(q\tau_{2})}{\Arcsin \tau_{2}}
$$
or, with $ x_{1} := \Arcsin \tau_{1} $, as
$$
\frac{\arcsin(q\Sin x_{1})}{x_{1}} <
\frac{\arcsin(q\Sin x_{2})}{x_{2}}.
$$
Since $ 0 < x_{1} < x_{2} $ this is true by Lemma 
\ref{lem_polyrig_dch}.
\QED

\medskip
\BL
\label{lem_polyrig_duniq}
Solely by the sidelengths of an oriented-convex
$ n $-gon whose vertices lie on a distance circle, the radius 
of the circle and the type of the polygon (niche or full
position) is determined.
\EL

\Bew
Consider two oriented-convex $ n $-gons with vertices on
distance circles of radii $ R, R' $ and with same
sidelengths $ L_{1},\ldots,L_{n} $. One has to identify 
$ R, R' $ and the two position types. The discussion follows
the eventual combinations of the positions.
\begin{DES1}{(iii)}
\item[(i)]
\TI{Circum-radii for two oriented-convex $ n $-gons with
same sidelengths in full position:}
This case is very easy. Namely, a division of 
\eqref{eq_polyrig_dvbil2} by $ 2\Sinh R  $ yields
\BE
\label{eq_ext_v1}
\pi = 
\sum_{k=1}^{n} 
\Arcsin\left(\frac{\ds\Sinh \frac{L_{k}}{2}}{\Sinh R}\right),
\EE
and this equation cannot allow two solutions for $ R $
because of the strict monotony of all summands w.r.t. $ R $.
\item[(ii)]
\TI{Circum-radii for two oriented-convex $ n $-gons with
same sidelengths in niche position:}
Assume again $ L_{n} > L_{k} $ for $ k = 1,\ldots,n-1 $.
From Corollary \ref{cor_polyrig_dmon} one deduces
\BE
\label{eq_ext_jn}
\frac{\Phi(R,L_{k})}{\Phi(R',L_{k})} > 
\frac{\Phi(R,L_{n})}{\Phi(R',L_{n})} := p,
\EE
If there existed two such polygons with circum-radii 
$ R' < R $ then, besides \eqref{eq_polyrig_dnbil1}, one had 
the analogue with $ R $ replaced by $ R' $:
\BE
\label{eq_polyrig_nischnisch}
\begin{aligned}
\Phi(R,L_{n}) &= \Phi(R,L_{1})+\cdots+\Phi(R,L_{n-1})\\[0.5ex]
\Phi(R',L_{n}) &= \Phi(R',L_{1})+\cdots+\Phi(R',L_{n-1}).
\end{aligned}
\EE
This would imply
$$
\frac{\Phi(R,L_{n})}{\Phi(R',L_{n})} =
\frac{\Phi(R,L_{1})+\cdots+\Phi(R,L_{n-1})}%
{\Phi(R',L_{1})+\cdots+\Phi(R',L_{n-1})} <
\frac{p\Phi(R',L_{1})+\cdots+p\Phi(R',L_{n-1})}%
{\Phi(R',L_{1})+\cdots+\Phi(R',L_{n-1})} = p,
$$
contradicting \eqref{eq_ext_jn}.
\item[(iii)]
\TI{Circum-radii for two oriented-convex $ n $-gons with
same sidelengths in the first mixed position:}
By this is meant that the polygon with circum-radius $ R $ is in
the full position and the polygon with circum-radius 
$ R' \leq R $ in the niche position. 

Assume $ R' < R $.
Then, on one hand, 
\eqref{eq_polyrig_dvbil2} holds and on the other hand 
\eqref{eq_polyrig_dnbil1} with $ R $ replaced by $ R' $:
$$
\begin{aligned}
2\pi \Sinh R &= \Phi(R,L_{1})+\cdots+\Phi(R,L_{n}) \\[0.5ex]
\Phi(R',L_{n}) &= \Phi(R',L_{1})+\cdots+\Phi(R',L_{n-1}).
\end{aligned} 
$$
With \eqref{eq_polyrig_dquot} this implies
$$
2\pi \Sinh R-\Phi(R,L_{n}) = 
\Phi(R,L_{1})+\cdots+\Phi(R,L_{n-1}) <
\Phi(R',L_{1})+\cdots+\Phi(R',L_{n-1}) =
\Phi(R',L_{n}),
$$
so
$$ 
2\pi \Sinh R' < 2\pi \Sinh R < 
\Phi(R,L_{n})+\Phi(R',L_{n}) < 
2\Phi(R',L_{n}),
$$
hence $ \pi \Sinh R' < \Phi(R',L_{n}) $. But this is
impossible by \eqref{eq_hypmot_LS2}: 
The first mixed position is at most possible for $ R = R' $. 
\item[(iv)]
\TI{Circum-radii for two oriented-convex $ n $-gons with
same sidelengths in the second mixed position:}
By this is meant that the polygon with circum-radius $ R $ is in
the niche position and the polygon with circum-radius 
$ R' \leq R $ in the full position.

Assume $ R' < R $.
Then, on one hand \eqref{eq_polyrig_dnbil1} holds and on the
other hand \eqref{eq_polyrig_dvbil2} with $ R $ replaced by 
$ R' $:
$$
\begin{aligned}
\Phi(R,L_{n}) &= \Phi(R,L_{1})+\cdots+\Phi(R,L_{n-1}) \\[0.5ex]
2\pi \Sinh R' &= \Phi(R',L_{1})+\cdots+\Phi(R',L_{n}).
\end{aligned}
$$
Since again $ L_{n} > L_{k} $ may be assumed for 
$ k = 1,\ldots,n-1 $ this implies by \eqref{eq_ext_jn}:
$$
\begin{aligned}
2\pi \Sinh R'-\Phi(R',L_{n}) &= 
\Phi(R',L_{1})+\cdots+\Phi(R',L_{n-1}) \\[0.5ex]
&<
\frac{1}{p}\Phi(R,L_{1})+\cdots+\frac{1}{p}\Phi(R,L_{n-1}) =
\frac{1}{p}\Phi(R,L_{n}) = \Phi(R',L_{n}),
\end{aligned}
$$
hence
$
2\pi \Sinh R' < 2\Phi(R',L_{n})
$,
what is impossible by \eqref{eq_hypmot_LS2}:
The second mixed position is at most possible for $ R = R' $. 
\end{DES1}

In fact, even for equal circum-radii $ R = R' $ the cases (iii) and
(iv) are not existent. If it were so, then consider a biggest
sidelength, say $ L_{n} $. Due to the niche position the
biggest sidelength is unique. Then both conditions 
\eqref{eq_polyrig_dnbil1} and \eqref{eq_polyrig_dvbil2} are 
satisfied simultaneously:
$$
\begin{aligned}
\Phi(R,L_{n}) &= \Phi(R,L_{1})+\cdots+\Phi(R,L_{n-1}) \\[0.5ex]
2\pi \Sinh R &= \Phi(R,L_{1})+\cdots+\Phi(R,L_{n}).
\end{aligned}
$$
This implies $ 2\pi\Sinh R = 2\Phi(R,L_{n}) $, so 
$ L_{n} = 2R $. The corresponding edge is then a diameter of
the distance circle, so contains the center $ O $, in
contrast to the niche position.
\QED

\medskip
In order to ensure the uniqueness of the laying off procedure in the 
case of a circumscribed distance circle, the following
statements are needed:

\medskip
\BL
\label{lem_polyrig_kreisp}
~\\[-3ex]
\begin{DES1}{(ii)}
\item[(i)]
For every point  $ Z \in \Sc(R) $ and every number 
$ L \in \left]0,2R\right[ $ there exists exactly one point 
$ W \in \Sc(R) $ with $ d(Z,W) = L $ such that the center 
$ O $ is situated left, resp. right of the directed line 
$ \overrightarrow{Z\vee W} $.
\item[(ii)]
Let $ Z,W,Z',W' $ be points in $ \Sc(R) $ with 
$ 0 < d(Z,W) = d(Z',W') < 2R $. If the center $ O $ is
situated either left of the directed lines 
$ \overrightarrow{Z\vee W} $ and 
$ \overrightarrow{Z'\vee W'} $ or right of them, then there 
exists a proper rotation $ D $ of $ \B $ around 
$ O $ with  $ D(Z) = D(Z') $ and $ D(W) = D(W') $.
\end{DES1}
\EL

\Bew
Of course, distances and the left and right positions of 
$ O $ are invariant under proper rotations. We only discuss the
case of the left position of $ O $. The other case runs 
analogously.

\smallskip
\TI{For (i):}
On account of this invariance one may assume, by applying of
a proper rotation, that 
$ z = \gamma(R,0) $ ($ \gamma $ as in \eqref{eq_hypmot_LS1}).
The conditions for $ w = \gamma(R,\vhi) $ with 
$ \vhi \in \left]-\pi,\pi\right[ $ then sound equivalently:
$$
\begin{aligned}
d(Z,W) = L & \iff
\Cos\vhi = 
1-2\left(\frac{\ds\Sinh \frac{L}{2}}{\Sinh R}\right)^{2}
\\
[\gamma(R,0),\gamma(R,\vhi),0] > 0 &\iff \Sin \vhi > 0.
\end{aligned}
$$
The first line follows from 
$ \norm{\gamma(R,\vhi)} = 1 $ and 
$ \Cosh d(z,w) = \cosh^{2}R-\sinh^{2}R\cdot\Cos\vhi $ with
means of identities for the hyperbolic functions. 

These conditions are satisfied by exactly one $ \vhi $.
Hence the uniqueness.

\smallskip
\TI{For (ii):}
By applying another proper rotation to $ Z' $ one may assume
$ z = z' = \gamma(R,0) $. The conditions on 
$ w = \gamma(R,\vhi) $ and 
$ w' = \gamma(R,\vhi' ) $ with 
$ \vhi,\vhi' \in\left]-\pi,\pi\right[ $ then sound: 
$$
\begin{aligned}
d(Z,W) = d(Z',W') &\iff 
\Cos \vhi = \Cos \vhi'\\
[\gamma(R,0),\gamma(R,\vhi),0] > 0, \; 
[\gamma(R,0),\gamma(R,\vhi'),0] > 0 &\iff 
\Sin \vhi > 0, \; \Sin \vhi' > 0.
\end{aligned}
$$
This implies $ \vhi = \vhi' $, so $ W = W' $, hence 
$ D := \id $ does it.
\QED

\medskip
Now the distance circle situation can be finished:

\medskip
\BT
\label{thm_polyrig_dkongr}
Any two oriented-convex $ n $-gons with vertices on 
distance circles and with the same length spectrum 
$ L_{1},\ldots,L_{n} $ are hyperbolically equivalent.
\ET

\Bew
By Lemma \ref{lem_polyrig_duniq} one can assume that both
polygons $ Z_{1}\ldots Z_{n} $ and $ Z_{1}'\ldots Z_{n}' $
have the same circum-distance-circle $ \Sc(R) $ and
that they are of the same position type.

If both polygons are in full position then the edgelines
$ \overrightarrow{Z_{1}\vee Z_{2}} $ and 
$ \overrightarrow{Z_{1}' \vee Z_{2}'} $ have the origin on
their left hand sides. So, by  Lemma  
\ref{lem_polyrig_kreisp} (ii), an additional rotation
may be applied on the second polygon in order to reach 
$ Z_{1} = Z_{1}' $ und $ Z_{2} = Z_{2}' $. It remains to
show $ Z_{k} = Z_{k}' $ for $ k = 3,\ldots,n $:
The points $ Z_{3} $ und $ Z_{3}' $ have the distance 
$ L_{2} \leq 2R $ from $ Z_{2} $, and the edgelines
$ \overrightarrow{Z_{2}\vee Z_{3}} $ and 
$ \overrightarrow{Z_{2}'\vee Z_{3}'} $ have the center 
$ O $ on their left hand side. In case $ L_{2} < 2R $
follows $ Z_{3} = Z_{3}' $ by Lemma 
\ref{lem_polyrig_kreisp} (i). If $ L_{2} = 2R $ then 
$ Z_{3} $ and $ Z_{3}' $ are the antipode of $ Z_{2} $,
hence equal. Obviously, one can proceed inductively in this
manner, finally obtaining the incidence of all vertices.

If both polygons are in niche position and if $ L_{n} $ is
the biggest sidelength (belonging to the edgelines
$ \overrightarrow{Z_{n}\vee Z_{1}} $ and
$ \overrightarrow{Z_{n}'\vee Z_{1}'} $) then these edgelines
have the center $ O $ on their right hand side. So, by
Lemma \ref{lem_polyrig_kreisp}(ii) one can achieve:
$ Z_{n} = Z_{n}' $ und $ Z_{1} = Z_{1}' $. Then 
$ Z_{2} $ and  $ Z_{2}' $ have the distance $ L_{1} < 2R $
from $ Z_{1} $, and the edgelines 
$ \overrightarrow{Z_{1}\vee Z_{2}} $ and 
$ \overrightarrow{Z_{1}'\vee Z_{2}'} $ have $ O $ on their
left hand side. By Lemma \ref{lem_polyrig_kreisp}(i) one
deduces $ Z_{2} = Z_{2}' $. Again, this reasoning can be
iterated until all vertices are identified.
\QED

\medskip
Now all components for the general statement are collected 
together:

\medskip
\BT[rigidity]
\label{thm_polyrig_haupt}
Two oriented-convex and cocyclic $ n $-gons of the
hyperbolic plane are hyperbolically equivalent if and only 
if they have the same length spectrum 
$ L_{1},\ldots,L_{n} $.
\ET

\Bew
Of course, the hyperbolic equivalence implies the same
sidelengths.
The converse follows from the synopsis of the Theorems
\ref{thm_polyrig_utyp}, \ref{thm_polyrig_gkongr},
\ref{thm_polyrig_akongr}, \ref{thm_polyrig_dkongr}.
\QED

\medskip
The oriented-convex and cocyclic polygons are the only
non-collinear polygons for which such a rigidity statement
can be hoped for. Namely, it will be seen in part 2 that for
any other non-collinear polygon there is an oriented-convex
and cocyclic polygon with same sidelengths but bigger area.
So the former polygon cannot be rigid. For collinear
polygons the rigidity can be discussed within the context of
the next section.


\section{\hspace{-1em}.
Converse of the $ n $-inequalities}
\label{polyex}
\markright{\ref{polyex}. Converse of the $ n $-inequalities}

The main theorem on rigidity is a result on uniqueness. In
its formulation above it doesn't contribute to the existence
of polygons with given sidelengths, though in the course of
the proof there were contained pieces of existence in form
of the lay off procedure. In fact, these will again enter the
game here. But the decisive mean will be the generalized
triangle inequality. Its behaviour will control the existence
question and also the rigidity for collinear polygons. 
The \TI{$ n $-inequalities} for positive real numbers 
$ L_{1},\ldots,L_{n} $ say:
\BE
\label{eq_polyex_inv1}
L_{k} \leq L_{1}+\cdots + \Hat{L_{k}} + \cdots + L_{n} 
\quad \Forall k = 1,\ldots,n.
\EE
If in all these inequalities appears the less than sign we
are speaking of the \TI{strict} case, otherwise of the 
\TI{non-strict} case.
Of course, \eqref{eq_polyex_inv1} is necessarily satisfied if
the $ L_{k} $ are the sidelengths of a polygon. In order to 
discuss the converse one first remarks that these
inequalities can be contracted to one inequality if the
numbering is such that $ L_{n} $ is the maximum of the 
$ L_{k} $. Then \eqref{eq_polyex_inv1} is in fact equivalent
to 
\BE
\label{eq_polyex_vor1}
L_{n} \leq L_{1}+\cdots+L_{n-1}. 
\EE

\TI{First, the strict case is discussed.}
(The non-strict case will be dealt with in the proof of
Theorem \ref{thm_polyex_inv}.)

The main means for the existence are the closedness
conditions from Sect. \ref{rigi}:
\begin{align}
\label{eq_polyex_G}
\tag{H}
\Sinh\frac{L_{n}}{2} 
&=
\Sinh\frac{L_{1}}{2}+\cdots+\Sinh\frac{L_{n-1}}{2}\\[1ex]
\label{eq_polyex_A}
\tag{D}
\Arsinh
\left(
\frac{\ds \Sinh\frac{L_{n}}{2}}{\Cosh\delta}
\right)
&=
\Arsinh
\left(
\frac{\ds \Sinh\frac{L_{1}}{2}}{\Cosh\delta}
\right)
+\cdots+
\Arsinh
\left(
\frac{\ds \Sinh\frac{L_{n-1}}{2}}{\Cosh\delta}
\right)\\[1ex]
\label{eq_polyex_DN}
\tag{CN}
\Arcsin
\left(
\frac{\ds \Sinh\frac{L_{n}}{2}}{\Sinh R}
\right)
&=
\Arcsin
\left(
\frac{\ds \Sinh\frac{L_{1}}{2}}{\Sinh R}
\right)
+\cdots+
\Arcsin
\left(
\frac{\ds \Sinh\frac{L_{n-1}}{2}}{\Sinh R}
\right)\\[1ex]
\label{eq_polyex_DV}
\tag{CF}
\pi 
&=
\Arcsin
\left(
\frac{\ds \Sinh\frac{L_{1}}{2}}{\Sinh R}
\right)
+\cdots+
\Arcsin
\left(
\frac{\ds \Sinh\frac{L_{n}}{2}}{\Sinh R}
\right).
\end{align}
See \eqref{eq_polyrig_agbil2}, \eqref{eq_polyrig_agbil1},
\eqref{eq_polyrig_dnbil1}, \eqref{eq_polyrig_dvbil2} und
\eqref{eq_hypmot_LS8}, \eqref{eq_hypmot_LS5}, 
\eqref{eq_hypmot_LS2}.

Given the lengths $ L_{1},\ldots,L_{n} $ as positive real
numbers, the last three equations are to be viewed as
conditions for the unknowns $ \delta > 0 $ and $ R > 0 $.
Since only certain combinations of quantities occur in these
equations we may set:
\BE
\label{eq_polyex_unk}
x := \frac{1}{\Cosh \delta} \quad 
\text{for equation \eqref{eq_polyex_A}},
\qquad
x := \frac{1}{\Sinh R} \quad 
\text{for equations \eqref{eq_polyex_DN} and
\eqref{eq_polyex_DV}}
\EE
and, generally
\BE
\label{eq_polyex_defl}
\lambda_{k} := \Sinh\frac{L_{k}}{2}.
\EE
Eqns. \eqref{eq_polyex_G}, \eqref{eq_polyex_A}, 
\eqref{eq_polyex_DN}, \eqref{eq_polyex_DV} are thus modified
to
\begin{align}
\label{eq_polyex_g}
\tag{h}
\lambda_{n}
&=
\lambda_{1}+\cdots+\lambda_{n-1} &&\\[1.5ex]
\label{eq_polyex_a}
\tag{d}
\arsinh(\lambda_{n}x) 
&=
\arsinh(\lambda_{1}x) +\cdots+\arsinh(\lambda_{n-1}x),
& 0 &< x < 1\\[1.5ex]
\label{eq_polyex_dn}
\tag{cn}
\arcsin(\lambda_{n}x) 
&=
\arcsin(\lambda_{1}x) +\cdots+\arcsin(\lambda_{n-1}x)),
& 0 &< x < \frac{1}{\lambda_{n}}\\[1ex]
\label{eq_polyex_dv}
\tag{cf}
\pi 
&=
\arcsin(\lambda_{1}x) +\cdots+\arcsin(\lambda_{n}x),
& 0 &< x \leq \frac{1}{\lambda_{n}}.
\end{align}
It must be proved that, for given positive reals
$ L_{1},\ldots,L_{n} $ with
\BE
\label{eq_polyex_vorst}
L_{j} \leq L_{n} \quad\Forall j, \qquad
L_{n} < L_{1}+\cdots+L_{n-1}. 
\EE
at least one of these four equations
has a solution $ x $ in the corresponding interval (the first
equation doesn't contain any unknown).

The last two equations have the same form as in the
Euclidean case, however with a different meaning for the
\TI{given} quantities $ \lambda_{1},\ldots,\lambda_{n} $. On their
solvability the following is known:

\medskip
\BL[$\text{Pinelis [2005]}$]
\label{lem_polyex_incr}
Let $ \lambda_{1},\ldots,\lambda_{n} \in \R $ be given,
satisfying $ 0 < \lambda_{j} \leq \lambda_{n} $ for
$ j = 1,\ldots,n $ and 
$ \lambda_{n} < \lambda_{1}+\cdots +\lambda_{n-1} $. Then:
\begin{DES1}{(ii)}
\item[(i)]
The function 
$$
g(x) := 
\arcsin(\lambda_{1}x)+\cdots+\arcsin(\lambda_{n-1}x)
-\arcsin(\lambda_{n}x)
$$
is positive in a certain interval right of $ 0 $. 
\item[(ii)]
If Eqn. \eqref{eq_polyex_dn} has \underline{no} solution in 
$ \left]0,1/\lambda_{n}\right[ $ then Eqn. 
\eqref{eq_polyex_dv} \underline{has} a solution in 
$ \left]0,1/\lambda_{n}\right] $.
\end{DES1}
\EL

Slightly modified, this can be shown as follows:

\medskip
\TI{Proof of Lemma \ref{lem_polyex_incr}.}

\TI{For (i):}
For positive $ x $ with 
$ \lambda_{1}x+\cdots+\lambda_{n-1}x \leq \frac{\pi}{2} $
by \eqref{rem_polyrig_umk2}:
$$
g(x) > 
\arcsin(\lambda_{1}x+\cdots+\lambda_{n-1}x)-
\arcsin(\lambda_{n}x).
$$
From $ \lambda_{1}x+\cdots+\lambda_{n-1}x > \lambda_{n}x $
follows $ g(x) > 0 $ for these $ x $. So the function 
$ g $ is positive in the interval $ \left]0,\beta\right[ $
where
$$
\beta:=
\frac{\pi}{2(\lambda_{1}+\cdots+\lambda_{n-1})}.
$$

\TI{For (ii):}
Since Eqn.  \eqref{eq_polyex_dn} has no solution in 
$ \left]0,1/\lambda_{n}\right[ $, the function 
$ g $ has a fixed sign in this interval, by part (i) a fixed
\TI{positive} sign. This implies
$$
g\left(\frac{1}{\lambda_{n}}\right) \geq 0.
$$
By means of $ g(x) $, Eqn. \eqref{eq_polyex_dv} sounds
equivalently
$$
g_{1}(x) := \pi-2\arcsin(\lambda_{n}x)-g(x) = 0.
$$
Now
$$ 
g_{1}(0) = \pi, \qquad 
g_{1}\left(\frac{1}{\lambda_{n}}\right) = 
-g\left(\frac{1}{\lambda_{n}}\right) \leq 0.
$$
Thus the equation $ g_{1}(x) = 0 $ has a solution in 
$ \left]0,1/\lambda_{n}\right] $.
\QED

\medskip
\BL
\label{lem_polyex_gadnv}
~\\[0.3ex]
Assume $ 0 < \lambda_{k} \leq \lambda_{n} $ for 
$ k = 1,\ldots,n $ and,
with the binding \eqref{eq_polyex_defl},
$ L_{n} < L_{1}+\cdots +L_{n-1} $.
Then at least one of the equations  \eqref{eq_polyex_g}, 
\eqref{eq_polyex_a},\eqref{eq_polyex_dn}, and 
\eqref{eq_polyex_dv} has a solution in the corresponding
interval.
\EL

\Bew
In case $ \lambda_{n} = \lambda_{1}+\cdots +\lambda_{n-1} $ 
it is Eqn. \eqref{eq_polyex_g}. In case 
$ \lambda_{n} < \lambda_{1}+\cdots +\lambda_{n-1} $, one of 
the equations \eqref{eq_polyex_dn}, \eqref{eq_polyex_dv} has
a solution, by Lemma \ref{lem_polyex_incr} (ii). Now, in case 
$ \lambda_{n} > \lambda_{1}+\cdots +\lambda_{n-1} $, it will
be shown that Eqn. \eqref{eq_polyex_a} is solvable: Define
$$
h(x) := 
\arsinh(\lambda_{n}x) 
-
\arsinh(\lambda_{1}x) -\cdots-\arsinh(\lambda_{n-1}x),
\quad x \in \R.
$$
One has $ h(0) := 0 $ and 
$ h'(0) = 
\lambda_{n}-\lambda_{1}-\cdots-\lambda_{n-1} > 0
$.
So $ h(x) $ is positive in a certain interval right of $ 0 $. 
On the other hand, by \eqref{eq_polyex_defl}
$$
h(1) = 
\frac{L_{n}}{2}-\frac{L_{1}}{2}-\cdots-\frac{L_{n-1}}{2}
< 0.
$$
Thus, $ h $ must vanish somewhere in $ \left]0,1\right[ $. 
\QED

\medskip
\medskip
\BT
\label{thm_polyex_inv}
Let $ L_{1},\ldots,L_{n} $ be positive real numbers with 
\BE
\label{eq_polyex_inv11}
L_{k} \leq L_{1}+\cdots + \Hat{L_{k}} + \cdots + L_{n} 
\quad \Forall k = 1,\ldots,n.
\EE
\vspace{-4ex}
\begin{DES1}{(ii)}
\item[(i)]
If in all inequalities \eqref{eq_polyex_inv11} occurs 
the strict less than symbol then there exists an
oriented-convex 
cocyclic $ n $-gon in the hyperbolic plane $ \B $ with
sidelengths $ L_{1},\ldots,L_{n} $. According to Theorems 
\ref{thm_polyrig_gkongr}, 
\ref{thm_polyrig_akongr}, and
\ref{thm_polyrig_dkongr}
such a polygon is uniquely determined up to hyperbolic
equivalence.
\item[(ii)]
If in one of the inequalities \eqref{eq_polyex_inv11} occurs
the equals sign, say in $ L_{n} = L_{1}+\cdots+L_{n-1} $,
then there also exists a $ n $-gon in the hyperbolic plane 
$ \B $ with sidelength $ L_{1},\ldots,L_{n} $. Such a
polygon is uniquely determined up to hyperbolic congruence, it
is collinear and its vertices $ Z_{1},\ldots,Z_{n} $ are in 
this order strictly monotonic arranged on the supporting
line, in particular pairwise distinct.
\end{DES1}
\ET

\Bew

\TI{For (i):}
Assume again $ L_{j} \leq L_{n} $ for $ j = 1,\ldots,n $ and
$ L_{n} < L_{1}+\cdots+L_{n-1} $.

The proof runs by inspecting the solvability of Eqns. 
\eqref{eq_polyex_g}, \eqref{eq_polyex_a},
\eqref{eq_polyex_dn}, \eqref{eq_polyex_dv} one by one. 
In all four cases it will be possible to obtain from a
solution $ x $ the corresponding parameter 
$ \delta $ resp. $ R $ via the substitutions
$$
\lambda_{k} := \Sinh\frac{L_{k}}{2}, \qquad
x = \frac{1}{\Cosh \delta}, \quad \text{resp.} \quad
x = \frac{1}{\Sinh R}. 
$$
Then, on the appropriate circle type, the lay off procedure
can be executed with the lengths $ L_{1},\ldots,L_{n-1} $ in
order to generate the vertices $ Z_{1},\ldots,Z_{n} $ of the
desired polygon. Thereby, the closedness condition ensures
each time the compatibility, i.e. the fact that also the
last distance $ d(Z_{n},Z_{1}) $ fits to $ L_{n} $.

\smallskip
\TI{Case i.h: Eqn. \eqref{eq_polyex_g} holds.}

This is the easiest case because no unknown $ x $ is
involved. Here, of course, the lay off is applied to the horocycle
$ \Hc $ as described in the proof of Theorem 
\ref{thm_polyrig_gkongr}: One parametrizes 
the orbit by arclength (positively constant proportional to 
the original group parameter $ b $) as $ \Gamma(s) $,
selects $ Z_{1} = \Gamma(s_{0}) $ arbitrarily and constructs
the further vertices by the rule
\BE
\label{eq_polyex_inv2}
Z_{k} :=
\Gamma(s_{0}+\Lambda(L_{1})+\cdots+\Lambda(L_{k-1})),
\quad k = 2,\ldots,n.
\EE
For $ k = n $ by Eqn. \eqref{eq_polyex_g} then holds
$$
\Lambda(L_{1})+\cdots+\Lambda(L_{n-1}) =
2\Sinh\frac{L_{1}}{2}+\cdots+2\Sinh\frac{L_{n-1}}{2} =
2\lambda_{1}+\cdots+2\lambda_{n-1} = 2\lambda_{n} =
\Lambda(L_{n}),
$$
so
\BE
\label{eq_polyex_inv3}
Z_{n} = \Gamma(s_{0}+\Lambda(L_{n})).
\EE
In case $ k \leq n-1 $ follows from \eqref{eq_polyex_inv2} 
for the the arclength between 
$ Z_{k} $ und $ Z_{k+1} $: $ s_{k} = \Lambda(L_{k}) $, hence by 
$ s_{k} = \Lambda(d(Z_{k},Z_{k+1})) $:
$ d(Z_{k},Z_{k+1}) = L_{k} $. In case $ k = n $ the same
results from Eqn. \eqref{eq_polyex_inv3}: 
$ d(Z_{1},Z_{n}) = L_{n} $.

As a result, the polygon $ Z_{1}\ldots Z_{n} $ fulfills all
requirements: it has $ \Hc $ as a circum-path, it is
oriented-convex by the strict monotony of the parameter
values, and the sidelengths are the given numbers
$ L_{1},\ldots,L_{n} $.

\smallskip
\TI{Case i.d: 
$ \lambda_{1}+\cdots+\lambda_{n-1} < \lambda_{n} $.}

So Eqn. \eqref{eq_polyex_a} has a solution $ x $ in the
given interval. The only difference to case i.h is that from
this solution a $ \delta > 0 $ must be determined by 
$$
\Cosh\delta = \frac{1}{x}.
$$
The lay off procedure now is to perform on the distance line
$ \Dc(\delta) $. The creation of the points 
$ Z_{1},\ldots,Z_{n} $ happens as above, of course with 
$ \Lambda $ replaced by $ \Psi $. By Eqn. 
\eqref{eq_polyex_a}, this time holds:
\BE
\label{eq_polyex_inv6}
\begin{aligned}
\Psi(\delta,L_{1})+\cdots+\Psi(\delta,L_{n-1}) &= 
2\Cosh\delta\Arsinh\frac{\lambda_{1}}{\Cosh\delta}
+\cdots+
2\Cosh\delta\Arsinh\frac{\lambda_{n-1}}{\Cosh\delta} \\
&=
2\Cosh\delta\Arsinh\frac{\lambda_{n}}{\Cosh\delta} =
\Psi(\delta,L_{n}),
\end{aligned}
\EE
with the additional arguments as above.

\smallskip
\TI{Case i.cn: 
$ \lambda_{1}+\cdots+\lambda_{n-1} > \lambda_{n} $
and Eqn. \eqref{eq_polyex_dn} is solvable.}

Here the arguing is somewhat different since arclength and 
chordlength are not always in a one to one relation. First
$ x $ is fixed as a solution of Eqn.  \eqref{eq_polyex_dn}
in the given interval and then $ R > 0 $ is determined by
$$
\Sinh R = \frac{1}{x}.
$$
The laying off is done on the distance circle $ \Sc(R) $.
Its arclength parametrization follows from Eqn. 
\eqref{eq_hypmot_LS1} (without changing the name) as
\BE
\label{eq_polyex_bpar}
\gamma(s) =
\BM
\Cosh R \\
\ds\Sinh R\cdot\Cos \frac{s}{\Sinh R} \\[2ex]
\ds\Sinh R\cdot\Sin \frac{s}{\Sinh R}
\EM,
\EE
and, for each arguments $ s_{1},s_{2} \in \R $, one
calculates from this
\BE
\label{eq_polyex_bdist}
\Sinh\frac{d(\Gamma(s_{1}),\Gamma(s_{2}))}{2} =
\Sinh R \cdot \Sin \frac{\abs{s_{2}-s_{1}}}{2\Sinh R}.
\EE
In order to puncture $ \Sc(R) $ suitably, the points $ Z_{k} $ are
constructed on $ \Sc(R) $ as follows, using a value $ s_{0} > 0 $ 
still to be chosen:
\BE
\label{eq_polyex_inv4}
Z_{k} :=
\Gamma(s_{0}+\Phi(R,L_{1})+\cdots+\Phi(R,L_{k-1})),
\quad k = 2,\ldots,n.
\EE
For $ k = n $ one has, analogous to Eqn. \eqref{eq_polyex_inv6},
this time on account of Eqn. \eqref{eq_polyex_dn}:
\BE
\label{eq_polyex_inv7}
\Phi(R,L_{1})+\cdots+\Phi(R,L_{n-1}) = \Phi(R,L_{n}),
\EE
so
\BE
\label{eq_polyex_inv5}
Z_{n} = \Gamma(s_{0}+\Phi(R,L_{n})).
\EE
Since always $ \Phi(R,L_{k}) \leq \pi\Sinh R $, the $ s_{0} $
can be chosen such that $ s_{0}+\Phi(R,L_{n}) < 2\pi\Sinh R $.
Then the points $ Z_{1},\ldots,Z_{n} $ form an
oriented-convex polygon with vertices on the distance circle 
$ \Sc(R) $,
punctured at $ \Gamma(0) $. By means of Eqn. \eqref{eq_polyex_bdist}
one calculates from \eqref{eq_polyex_inv4}, resp.
\eqref{eq_polyex_inv5}:
$$
d(Z_{k},Z_{k+1}) = L_{k}, \quad k = 1,\ldots,n-1 
\qquad \text{resp.} \qquad
d(Z_{1},Z_{n}) = L_{n}.
$$
So the polygon constructed in this way has all desired
properties.

\smallskip
\TI{Case i.cf: 
$ \lambda_{1}+\cdots+\lambda_{n-1} > \lambda_{n} $
and Eqn. \eqref{eq_polyex_dv} is solvable.}

The construction runs analogously to the foregoing case with
the only difference that Eqn. \eqref{eq_polyex_inv7} has to 
be replaced by the following equation which rests on 
\eqref{eq_polyex_dv}:
$$
\Phi(R,L_{1})+\cdots+\Phi(R,L_{n-1}) = 2\pi\Sinh R-\Phi(R,L_{n}).
$$
Then
$$
Z_{n} = \Gamma(s_{0}2+2\pi\Sinh R-\Phi(R,L_{n})), 
$$
where this time one can achieve 
$ s_{0}+2\pi\Sinh R-\Phi(R,L_{n}) < 2\pi\Sinh R $ by an
appropriate choice of $ s_{0} > 0 $. The points 
$ Z_{1},\ldots,Z_{n} $ then form an oriented-convex polygon 
with vertices on the distance circle $ \Sc(R) $,
punctured at $ \Gamma(0) $.

As to the distances, as in the foregoing case:
$ d(Z_{k},Z_{k+1}) = L_{k} $, $ k = 1,\ldots,n-1 $.
However, for $ k = n $ one has to argue differently. Again
by Eqn. \eqref{eq_polyex_bdist}:
$$
\begin{aligned}
\Sinh\frac{d(Z_{1},Z_{n})}{2} 
&=
\Sinh R \cdot \Sin \frac{\abs{2\pi\Sinh R-\Phi(R,L_{n})}}%
{2\Sinh R} \\[1ex]
&=
\Sinh R\cdot\sin
\left(\pi-\frac{\Phi(R,L_{n})}{2\Sinh R}\right) \\[1ex]
&=
\Sinh R\cdot\Sin\frac{\Phi(R,L_{n})}{2\Sinh R} \\
&=
\Sinh\frac{L_{n}}{2}.
\end{aligned}
$$
Also this time, the polygon $ Z_{1}\ldots Z_{n} $ has the desired 
properties.

\medskip
\TI{For (ii):}
Assume $ L_{n} = L_{1}+\cdots+L_{n-1} $. Essentially, the
assertions are well-known in this degenerate case.

\TI{Existence:}
The points 
$ Z_{1},\ldots, Z_{n} $ can be generated by
laying off successively the segments of lengths $ L_{1},\ldots,L_{n-1} $
along the groundline $ U_{0} $ in positive direction, starting e.g.
from $ Z_{1} := O $. The assumption then ensures that 
$ d(Z_{1},Z_{n}) = L_{n} $ because 
translations along $ U_{0} $ preserve the distance between
pre-image and image for points on $ U_{0} $.

\smallskip
\TI{Uniqueness:}
This amounts to show:

\TI{For points $ Z_{1},\ldots,Z_{n} \in \B $, the condition 
\BE
\label{eq_polyex_gln}
\tag{n}
d(Z_{1},Z_{n}) =
d(Z_{1},Z_{2})+d(Z_{2},Z_{3})+\cdots+d(Z_{n-1},Z_{n})
\EE
is necessary and sufficient for that $ Z_{1},\ldots,Z_{n} $ 
lie on a line and are monotonically arranged on it.}

The necessity is obvious. The sufficiency follows by
induction on $ n $. Without loss of generality, one may
assume $ Z_{k} \neq Z_{k+1} $ for $ k = 1,\ldots,n-1 $.
Otherwise identical successive points can be contracted to a 
single point. 

\TI{Initial step $ n = 3 $:}
From $ d(Z_{1},Z_{3}) = d(Z_{1},Z_{2})+d(Z_{2},Z_{3}) $
first follows $ Z_{1} \neq Z_{3} $. If $ Z_{2} $ were not on
$ Z_{1}\vee Z_{3} $ then the hyperbolic cosine rule would
imply $ d(Z_{1},Z_{3}) < d(Z_{1},Z_{2})+d(Z_{2},Z_{3}) $, a 
contradiction. So $ Z_{1},Z_{2},Z_{3} $ must be collinear.
Moreover, if $ Z_{2} $ were not on the segment 
then one had again 
$ d(Z_{1},Z_{3}) < d(Z_{1},Z_{2})+d(Z_{2},Z_{3}) $, so 
$ Z_{1},Z_{2},Z_{3} $ have to be collinear \TI{and} monotonically
arranged.

\TI{Induction step from $ n-1 $ to $ n $:}
First, one has
$ d(Z_{1},Z_{3}) = d(Z_{1},Z_{2})+d(Z_{2},Z_{3}) $. In fact,
if one had $ d(Z_{1},Z_{3}) < d(Z_{1},Z_{2})+d(Z_{2},Z_{3}) $,
then \eqref{eq_polyex_gln} would imply 
$$
d(Z_{1},Z_{n}) > 
d(Z_{1},Z_{3})+d(Z_{3},Z_{4})+\cdots+d(Z_{n-1},Z_{n}) 
\geq
d(Z_{1},Z_{n}),
$$
a contradiction. 
From 
$ d(Z_{1},Z_{3}) = d(Z_{1},Z_{2})+d(Z_{2},Z_{3}) $ follows
$ Z_{1} \neq Z_{3} $ and,
by the initial step, the monotonic arrangement of 
$ Z_{1},Z_{2},Z_{3} $ on a line and also
$$
d(Z_{1},Z_{n}) = 
d(Z_{1},Z_{3})+d(Z_{3},Z_{4})+\cdots+d(Z_{n-1},Z_{n}).
$$
By the induction hypothesis, this implies the monotonic
arrangement of $ Z_{1},Z_{3},Z_{4},\ldots,Z_{n} $ on a line.
Altogether, the monotonic arrangement of 
$ Z_{1},Z_{2},Z_{3},\ldots,Z_{n} $ on a line is thus
achieved.

In the present case of polygons, \eqref{eq_polyex_gln}
additionally implies the strictly monotonic arrangement
since all distances on the right hand side are positive.

If two such polygons $ Z_{1}Z_{2}\ldots Z_{n} $ and
$ Z_{1}'Z_{2}'\ldots Z_{n}' $ with same sidelengths 
$ L_{k} = d(Z_{k},Z_{k+1}) = d(Z_{k}',Z_{k+1}') $, 
$ k = 1,\ldots,n-1 $, are given then first the supporting lines
can be identified and then also the point 
$ Z_{1} $ with $ Z_{1}' $ and the point $ Z_{2} $ with 
$ Z_{2}' $, by suitable motions. By the equal sidelengths and 
the monotonic arrangement, then also the vertices 
$ Z_{k} $ will be identical with the corresponding $ Z_{k}' $
for $ k = 3,\ldots,n $.
\QED

\medskip
\BC[conversion of the generalized triangle inequality]
\label{cor_polyex_hn}
Given, for $ n \geq 3 $, positive real numbers 
$ L_{1},\ldots,L_{n} $ then the inequalities
\BE
\label{eq_polyex_hn1}
L_{k} \leq L_{1}+\cdots + \Hat{L_{k}} + \cdots + L_{n},
\qquad k = 1,\ldots,n,
\EE
are necessary and sufficient for that there exists a 
$ n $-gon in the hyperbolic space $ \H^{m} $ of dimension 
$ m \geq 2 $ with sidelengths $ L_{1},\ldots,L_{n} $.
\EC

\Bew

The \TI{necessity} is in order by the generalized triangle
inequality.

For the \TI{sufficiency}, one only has to observe that 
the hyperbolic space $ \H^{m} $ always contains one,
in fact many, hyperbolic planes. So, already in such a plane
there exists a polygon with the given sidelengths. It can be
chosen collinear resp. cocyclic and oriented-convex
(with a deliberate orientation of the plane) if in the
inequalities \eqref{eq_polyex_hn1} once resp. never occurs the 
equals sign.
\QED

\bigskip
\TB{\Large References}

\smallskip
Blumenthal, L. M [1970]:
Theory and applications of distance geometry:
Second Edition,
Chelsea Publ. Comp. Bronx, New York, i-xi and 1-347

Coxeter, H. S. M. [1968]:
Non-Euclidean geometry:
reprint 1968 of the fifth edition,
Univ. Toronto Press, i-xv and 1-309

Klingenberg, W. [1978]:
A course in differential geometry:
English version of the German edition from 1973,
Springer-Verlag New York-Berlin-Heidelberg,
i-vii and 1-178

Leichwei"s, K. [2003]:
On the addition of convex sets in the hyperbolic plane:
J. Geom. 78, 92-121

Leichwei"s, K. [2004]:
Support function and hyperbolic plane:
Manuscr. Math.. 114, 177-196

Leichwei"s, K. [2005]:
Curves of constant width in the non-Euclidean geometry:
Abh. Math. Semin. Univ. Hamb. 75, 257-284

Leichwei"s, K. [2008.a]:
Polar curves in the non-Euclidean geometry:
Result. Math. 52, 143-160

Leichwei"s, K. [2008.b]:
On Steiner's symmetrization of convex bodies in 
non-Euclidean geometry:
Result. Math. 52, 339-346

Lenz, H. [1967]:
Nichteuklidische Geomertrie:
Bibliographisches Institut Mannheim, 1-235

Menger, K. [1928]:
Untersuchungen "uber allgemeine Metrik:
Math. Ann. 100, 75-163

Pinelis, I. [2005]:
Cyclic polygons with given edge lengths:
Existence and uniqueness:
J. Geom. 82, 156-171

\medskip
Rolf Walter\\
Fakult"at f"ur Mathematik\\
Technische Universit"at Dortmund\\
Arbeitsgebiet Differentialgeometrie\\
Vogelpothsweg 87\\
D-44227 Dortmund

E-Mail: rolf.walter$ @ $tu-dortmund.de

\end{document}